\documentclass[3p,times]{elsarticle}

\usepackage{array}
\usepackage{color}
\usepackage{xcolor,colortbl}
\usepackage{natbib}
\usepackage{tabularx}
\usepackage{graphicx}
\usepackage{amsmath}
\usepackage{amssymb}
\usepackage{amsthm}
\usepackage{empheq}
\usepackage{multirow}
\usepackage{subfig}
\usepackage{moreverb}
\usepackage{nicefrac} 
\usepackage{dsfont}
\usepackage{bm}
\usepackage{multirow}
\usepackage{soul}

\newcommand\BibTeX{{\rmfamily B\kern-.05em \textsc{i\kern-.025em b}\kern-.08em
T\kern-.1667em\lower.7ex\hbox{E}\kern-.125emX}}

\usepackage{hyperref} 
\hypersetup{
    colorlinks=true,                          
    linkcolor=blue, 
    citecolor=red, 
    urlcolor=blue  } 

\newcommand{\Q}{\mathbf{Q}}
\renewcommand{\S}{\mathbf{S}}

\newcommand{\eg}{\textit{e.g.}}

\newcommand{\A}{\mathbf{A}}

\newcommand{\be}{\begin{equation}}
\newcommand{\ee}{\end{equation}}
\newcommand{\bdm}{\begin{displaymath}}
\newcommand{\edm}{\end{displaymath}}

\newcommand{\bea}{\begin{eqnarray} }
\newcommand{\eea}{\end{eqnarray} }

\definecolor{darkgreen}{rgb}{0.0, 0.5, 0.0}

\newfont{\numerikEleven}{ecrm1000}
\newfont{\numerikTen}{cmss10}
\newfont{\numerikNine}{cmss9}
\newfont{\numerikEight}{cmss8}

\journal{Journal of Computational Physics}

\begin{document} 
\begin{frontmatter} 
\title{On GLM curl cleaning for a first order reduction of the CCZ4 formulation of the Einstein field equations} 

\author[UniTN]{Michael Dumbser$^{*}$}
\ead{michael.dumbser@unitn.it}
\cortext[cor1]{Corresponding author}

\author[MPI]{Francesco Fambri}
\ead{francesco.fambri@ipp.mpg.de}

\author[UniTN]{Elena Gaburro}
\ead{elena.gaburro@unitn.it}

\author[TUM]{Anne Reinarz}
\ead{reinarz.in@tum.de}

\address[UniTN]{Department of Civil, Environmental and Mechanical Engineering, University of Trento, Via Mesiano 77, I-38123 Trento, Italy}
\address[MPI]{Max-Planck Institute for Plasma Physics, Boltzmannstr. 2, D-85748 Garching, Germany}
\address[TUM]{Department of Informatics, Technical University Munich (TUM), Boltzmannstr. 3, D-85748 Garching, Germany}

\begin{abstract}

In this paper we propose an extension of the generalized Lagrangian multiplier method (GLM) of Munz et al.  \cite{MunzCleaning,Dedneretal}, which was originally conceived for the numerical solution of the Maxwell and MHD equations with divergence-type involutions, to the case of hyperbolic PDE systems with curl-type involutions. {The key idea here is to solve an \textit{augmented} PDE system, in which curl errors propagate away via a Maxwell-type evolution system. The new approach is first presented on a simple model problem, in order to explain the basic ideas. Subsequently, } we apply it to a strongly hyperbolic first order reduction of the CCZ4 formulation (FO-CCZ4) 
of the Einstein field equations of general relativity, which is endowed with 11 curl constraints. 
Several numerical examples, including the long-time evolution of a stable neutron star in anti-Cowling approximation, 
are presented in order to show the obtained improvements with respect to the standard formulation without special treatment of the curl involution constraints. 

{The main advantages of the proposed GLM approach are its complete independence of the underlying numerical scheme and grid topology and its easy implementation into existing computer codes. However, this flexibility comes at the price of needing to add for each curl involution one additional 3 vector plus another scalar in the augmented system for homogeneous curl constraints, and even two additional scalars for non-homogeneous curl involutions. For the FO-CCZ4 system with 11 homogeneous curl involutions, this means that additional 44 evolution quantities need to be added. }

\end{abstract}

\end{frontmatter}


%
\section{Introduction} \label{sec.introduction}

Hyperbolic PDE systems with involutions are very frequent in science and engineering. The most common involutions are linear and of the divergence and curl type, the most famous involution being the divergence-free condition of the magnetic field in the Maxwell and the magnetohydrodynamics (MHD) equations.  
While a lot of research has been dedicated in the past to the appropriate numerical discretization of PDE with divergence constraints, much less is known on curl-preserving numerical schemes for PDE with curl involutions. 
However, nowadays there exists an increasing number of hyperbolic PDE systems in mechanics and physics that is endowed with natural curl involutions, where the curl of a vector field should remain zero for all times, if it was initially zero.
The most prominent examples are the system of nonlinear hyperelasticity of Godunov, Peshkov and Romenski (GPR model) \cite{GodunovRomenski72,PeshRom2014,GPRmodel} written in terms of the distortion field $\mathbf{A}$ and the conservative compressible multi-phase flow model of Romenski et al. \cite{Rom1998,RomenskiTwoPhase2010}. All the aforementioned mathematical models fall into the larger class of symmetric hyperbolic and thermodynamically compatible (SHTC) systems, studied by Godunov and  Romenski et al. in \cite{Godunov1961,Rom1998,Godunov:2003a,SHTC-GENERIC-CMAT}.  
Further systems with curl involutions include the new 
hyperbolic model for surface tension and the recent hyperbolic reformulation of the Schr\"odinger  equation of Gavrilyuk and Favrie et al.  \cite{Schmidmayer2016,Dhaouadi2018}, as well as first order reductions of the Einstein field equations such as those proposed, e.g., in \cite{Alic:2009,Brown2012,ADERCCZ4}. Finally, we also would like to point out that the model of Newtonian continuum physics including solid mechanics and electro-dynamics proposed and discretized in \cite{Rom1998,GPRmodelMHD} contains both, curl and divergence involutions. Furthermore, the Einstein field equations in 3+1 ADM split \cite{Arnowitt59} include also \textit{nonlinear} second order involutions for the four-dimensional metric tensor, which are the well-known ADM constraints, namely the Hamiltonian constraint and the three momentum constraints, see \cite{Alcubierre:2008}. These nonlinear involutions are typically treated either by adding suitable multiples
of the constraint to the governing PDE system, similar to the Godunov-Powell term in MHD \cite{God1972MHD,PowellMHD1,PowellMHD2}, 
or via the so-called Z4 constraint cleaning \cite{Alic:2009,Alic:2012}, which can be seen as a generalization of the GLM approach of Munz et al. \cite{MunzCleaning,Dedneretal} to the preservation of nonlinear involutions within the Einstein field equations.  

Here we briefly recall the hyperbolic GLM approach of Munz et al. \cite{MunzCleaning,Dedneretal} for the Maxwell and MHD equations. Throughout this paper we employ the Einstein summation convention, which implies summation over two repeated indices. Furthermore, we will use greek indices that range from 0 to 3 and latin indices ranging from 1 to 3. We furthermore use the notation 
$\partial_t = \partial / \partial t$,  $\partial_k = \partial / \partial x_k$ and $\varepsilon_{ijk}$ is the
usual fully anti-symmetric Levi-Civita symbol. 
The induction equation in electrodynamics reads 
\begin{equation}
\partial_t B_k + \varepsilon_{kij} \partial_i E_j = 0,  
\label{eqn.induction} 
\end{equation} 
with the magnetic field $B_k$ and the electric field $E_j$. A consequence of the above induction equation is the 
famous divergence-free condition 
\begin{equation}
\mathcal{I} = \partial_m B_m = 0
\label{eqn.divBzero}
\end{equation} 
of the magnetic field, which states that there exist no magnetic monopoles, or, in other words, the magnetic field
will remain divergence-free for all times if it was initially divergence-free. One way to preserve a divergence-free
magnetic field within a numerical scheme is the use of an exactly divergence-free discretization on appropriately  staggered meshes, see
e.g. \cite{Yee66,DeVore,BalsaraSpicer1999,Balsara2004,GardinerStone,balsarahllc2d,MUSIC1,ADERdivB}. However, the implementation of such exactly structure-preserving methods into an existing code is rather invasive and often requires substantial changes in the algorithm structure of  existing general purpose solvers for hyperbolic conservation laws that were not right from the beginning designed for the solution hyperbolic PDE with involution constraints.  
The very popular GLM method proposed by Munz et al. in \cite{MunzCleaning,Dedneretal} is an alternative to exactly constraint--preserving schemes and requires only a rather small change 
\textit{at the PDE level}, where simply an additional equation for a cleaning scalar $\phi$ is added to the system, so that divergence errors cannot accumulate locally any more, but instead are transported away under the form of acoustic-like waves with finite speed. This approach is very easy to implement in any general purpose CFD code and is completely independent of the underlying numerical scheme or mesh topology. The role of the cleaning scalar $\phi$ is the one of a generalized Lagrangian multiplier (GLM) that accounts for the involution constraint. The way how it works can best be explained with a physical example, which is the role of the pressure in the compressible Euler equations: for low Mach numbers (i.e. for large sound speed compared to the flow speed), the coupling of the momentum equation with the pressure equation drives the divergence of the velocity field to zero for $M \to 0$. In the same manner, the evolution equation of the additional cleaning scalar $\phi$ coupled with the induction equation drives the divergence of the magnetic field to zero if the cleaning speed is chosen large enough. The \textit{augmented induction equation} according to the GLM approach of Munz et al. \cite{MunzCleaning,Dedneretal} therefore reads 
\begin{eqnarray}
\label{eqn.induction.glm} 
\partial_t B_k + \varepsilon_{kij} \partial_i E_j + \textcolor{red}{ \partial_k \varphi } &=& 0, \\   
\label{eqn.phi.glm} 
\textcolor{red}{\partial_t \varphi + a_d^2  \, \partial_m B_m }  &=& \textcolor{red}{-\epsilon_d \varphi},  
\end{eqnarray} 
with the new cleaning scalar $\varphi$, an \textit{artificial} cleaning speed $a_d$ and a small damping parameter $\epsilon_d$. 
The new terms in the augmented system \eqref{eqn.induction.glm} and \eqref{eqn.phi.glm} with respect to the original
equation \eqref{eqn.induction} are highlighted in red, for convenience. It is easy to see that for $a_d \to \infty$
the equation \eqref{eqn.phi.glm} leads to $\partial_m B_m \to 0$, which is the above involution \eqref{eqn.divBzero}. 
The rest of this paper is structured as follows: in order to show the main idea on a simple and clear example, in Section \ref{sec.glm} we first present the extended GLM approach for hyperbolic PDE with curl-type involutions on a simple toy model. Next, in Section \ref{sec.model} we introduce the extended GLM approach for the full first order reduction of the 
CCZ4 formulation of the Einstein field equations (FO-CCZ4) with matter source terms. In Section \ref{sec.results} we present some numerical results that clearly show the benefits of the extended GLM method and in Section \ref{sec.conclusion} we give some concluding remarks and an outlook to future work.  

{At this point we clearly emphasize that in this paper we do \textit{not yet} consider the fully coupled Einstein-Euler system, i.e. the self-consistent evolution of the spacetime coupled with the matter via the general relativistic Euler or MHD equations. This is out of scope of the present manuscript and will be considered elsewhere. In this paper, we only consider the Einstein field equations in vacuum or with prescribed matter source terms, i.e. the so-called anti-Cowling approximation. }

\section{Hyperbolic curl cleaning with an extended generalized Lagrangian multiplier approach} 
\label{sec.glm} 

We first show the basic idea of our new approach on a simple toy model, in order to ease notation and to facilitate the understanding of the underlying concepts, before applying the method to the full Einstein field equations. Consider the following evolution system for {one scalar $\rho$ and } two vector fields $v_k$ and $J_k$,   
\begin{eqnarray}
    \label{eqn.rho}
\partial_t \rho + \partial_i \left( \rho v_i  \right) &=& 0,  \\ 
    \label{eqn.vk}
    \partial_t (\rho v_k) + \partial_i \left( \rho v_i v_k + \rho  c_0^2 J_i J_k \right) &=& 0,  \\ 
  \label{eqn.toy} 
  \partial_t J_k + \partial_k ( v_m J_m ) + v_m \left( \partial_m J_k - \partial_k J_m \right) &=& 0, 
\end{eqnarray} 
with $c_0$ a given constant. Defining a scalar quantity $\chi =  v_m J_m $ and with 
the use of the Schwarz theorem $\partial_k \partial_m \chi = \partial_m \partial_k \chi$ (symmetry of second derivatives) it is very easy to see that the second PDE above, eqn. \eqref{eqn.toy}, is endowed with the linear involution constraint  
\begin{equation}
   \mathcal{I}_{mk} = \partial_m J_k - \partial_k J_m = 0.  
\end{equation}
This means that if $\mathcal{I}_{mk}=0$ at the initial time, it will remain zero for all times. Indeed, one can
immediately notice that the involution itself is contained in the third term on the left hand side of 
eqn. \eqref{eqn.toy}. This term is \textit{very similar} to the so-called Godunov-Powell term in the MHD equations, see \cite{God1972MHD,PowellMHD1}, which was found by Godunov in 1972 in order to symmetrize the MHD system and which 
was later used by Powell in order to improve the behavior of numerical methods when discretizing the MHD equations in multiple space dimensions. 
For a general purpose numerical method applied to \eqref{eqn.toy}, it is very hard to guarantee $\mathcal{I}_{mk}=0$ at 
the discrete level. Satisfying the involution exactly would require a structure-preserving scheme, similar to those
employed for the Maxwell and MHD equations \cite{Yee66,DeVore,BalsaraSpicer1999} or those forwarded in 
\cite{HymanShashkov1997,JeltschTorrilhon2006,Torrilhon2004}. However, it might be very cumbersome 
to add such structure preserving schemes into an existing general purpose scheme for hyperbolic conservation laws, 
since these exactly involution satisfying methods usually require a particular staggering of the data on edges, 
faces and volumes and have thus to be foreseen right from the beginning when designing the scheme and the related software. 
Also, these structure preserving schemes might not be easy to implement on general meshes or 
for all types of high order discretizations (finite differences, ENO/WENO finite volume schemes, discontinuous Galerkin methods etc.). 

{Nowadays, exactly divergence--preserving discontinuous Galerkin and finite volume schemes are available at all orders on Cartesian grids, on structured curvilinear meshes, on unstructured simplex meshes and on geodesic meshes, see e.g.   
	\cite{BalsaraSpicer1999,BalsaraAMR,BalsaraMultiDRS,ADERdivB,balsarahlle2d,balsarahlle3d,MUSIC1,MUSIC2,BalsaraCED,BalsaraKaeppeli,BalsaraGeodesic,HazraBalsara,ShadabBalsara}. 
		However, much less is known so far on the construction of arbitrary high order accurate exactly  divergence-preserving schemes on general polygonal and polyhedral meshes \cite{Springel,ArepoTN,ShashkovMultiMat3}, or on general 
		space-trees with arbitrary refinement factor $\mathfrak{r}$, see e.g. \cite{AMR3DCL}.  }  
	
As already mentioned in the introduction, the main advantage of the GLM approach of Munz et al. 
\cite{MunzCleaning,Dedneretal} for divergence constraints is not only the ease of implementation, but also its great \textit{flexibility} and its compatibility with \textit{all} types of mesh topologies and numerical schemes, since it only requires the solution of 
an additional scalar PDE for the cleaning scalar, which can easily be added to an existing code.  

The extended GLM curl cleaning proposed in this paper can now be explained on the toy system  \eqref{eqn.vk}-\eqref{eqn.toy} as follows. The original governing PDE system \eqref{eqn.vk} and \eqref{eqn.toy} is 
simply \textit{replaced} by the following \textit{augmented system}   
\begin{eqnarray}
\label{eqn.toy.glm} 
\partial_t \rho + \partial_i \left( \rho v_i  \right) &=& 0,  \\ 
\partial_t (\rho v_k) + \partial_i \left( \rho v_i v_k + \rho c_0^2 J_i J_k \right) &=& 0,  \\ 
\label{eqn.glm.Jk} 
\partial_t J_k + \partial_k ( v_m J_m ) + v_m \left( \partial_m J_k - \partial_k J_m \right) + \textcolor{blue}{\varepsilon_{klm} \partial_l \psi_m} &=& 0, \\ 
\label{eqn.glm.psi} 
\textcolor{blue}{ 
	\partial_t \psi_k - a_c^2 \, \varepsilon_{klm} \partial_l J_m } \textcolor{red}{ + \partial_k \varphi } &=& 
\textcolor{blue}{  - \epsilon_{c} \, \psi_k,  }  \\  
\label{eqn.glm.phi} 
\textcolor{red}{ 
\partial_t \varphi + a_d^2 \, \partial_m \psi_m } &=& \textcolor{red}{- \epsilon_{d} \varphi}, 
\end{eqnarray}
where $a_c$ is a new cleaning speed associated with the curl cleaning. The new terms associated with the curl cleaning are highlighted in blue, for convenience, while the terms of the original PDE \eqref{eqn.toy} are written in black. Since the evolution equation for the cleaning vector field $\psi_k$ has formally the same structure as the induction equation \eqref{eqn.induction} of the Maxwell equations, it is again endowed with the divergence-free constraint 
$\partial_m \psi_m = 0$, which is taken into account via the classical GLM method (red terms).  
It is easy to see that from \eqref{eqn.glm.psi} for 
$a_c \to \infty$ we obtain $ \epsilon_{klm} \partial_l J_m \to 0$ in the limit, thus satisfying 
the involution in the sense $\mathcal{I}_{mk} \to 0$. The augmented system \eqref{eqn.toy.glm}-\eqref{eqn.glm.phi} can now be solved with any standard numerical method 
for nonlinear systems of hyperbolic partial differential equations. 

{In case the curl involutions are \textit{not} homogeneous, but where the curl has to assume a prescribed non-zero value, which is the case when the evolution equation \eqref{eqn.toy} for $J_k$ contains a source term $S_k$, 
\begin{equation} 
\partial_t J_k + \partial_k ( v_m J_m ) + v_m \left( \partial_m J_k - \partial_k J_m \right) = S_k,
\label{eqn.Jk.src} 
\end{equation} 		
	see \cite{Rom1998,Godunov:2003a,GPRmodel,GPRTorsion,SHTC-GENERIC-CMAT}, then the following strategy can be used: 
We first define the curl of $\mathbf{J}$ to be equal to another quantity 
\begin{equation}
  B_i = \varepsilon_{ijk} \partial_j J_k,  
  \label{burgers.involution} 
\end{equation} 
which according to \cite{SHTC-GENERIC-CMAT} is the so-called \textit{Burgers vector}. Its evolution equation can be obtained by  taking the curl of \eqref{eqn.Jk.src}, leading to the following additional evolution equation (see Appendix C of \cite{SHTC-GENERIC-CMAT}): 
\begin{equation}
\partial_t B_i + \partial_k \left( B_i v_k - v_i B_k - \varepsilon_{ikj} S_j \right) + v_i \, \partial_k B_k = 0. 
   \label{eqn.Burgers} 
\end{equation}
Directly from its definition \eqref{burgers.involution} it is obvious that the Burgers vector must be divergence-free, i.e. $\partial_k B_k = 0$. 
The divergence constraint on the Burgers vector is explicitly contained in \eqref{eqn.Burgers} in order to achieve a Galilean invariant formulation, see also \cite{God1972MHD,Rom1998,GPRmodelMHD,GPRTorsion}. 
Therefore, the GLM approach for a general, non-trivial curl of $J_k$ reads as follows:  }
\begin{eqnarray}
\partial_t \rho + \partial_i \left( \rho v_i  \right) &=& 0,  \\ 
\label{eqn.ex.v} 
\partial_t (\rho v_k) + \partial_i \left( \rho v_i v_k + \rho c_0^2 J_i J_k \right) &=& 0,  \\ 
\label{eqn.ex.J} 
\partial_t J_k + \partial_k ( v_m J_m ) + v_m \left( \partial_m J_k - \partial_k J_m \right) + \textcolor{blue}{\varepsilon_{klm} \partial_l \psi_m} &=& \textcolor{darkgreen}{S_k}, \\ 
\label{eqn.ex.B} 
\textcolor{darkgreen}{
   \partial_t B_i + \partial_k \left( B_i v_k - v_i B_k - \varepsilon_{ikj} S_j \right) + v_i \partial_k B_k}
+ \textcolor{red}{\partial_i \chi} &=&
\textcolor{darkgreen}{0},  \\  
\label{eqn.ex.psi} 
\textcolor{blue}{ 
	\partial_t \psi_k - a_c^2 \, \varepsilon_{klm} \partial_l J_m } \textcolor{red}{ + \partial_k \varphi } &=& 
\textcolor{blue}{  - \epsilon_{c} \, \psi_k  }  - \textcolor{darkgreen}{a_c^2 B_k} \\  
\label{eqn.ex.phi} 
\textcolor{red}{ 
	\partial_t \varphi + a_d^2 \, \partial_m \psi_m } &=& \textcolor{red}{- \epsilon_{d} \varphi}, \\
\textcolor{red}{
		\partial_t \chi + a_b^2 \, \partial_m B_m } &=& \textcolor{red}{- \epsilon_{b} \chi}. 
\end{eqnarray}
{
Here, we have used the usual blue color for the curl cleaning, the red color for the divergence cleaning and the green color due to the additional terms and equations that are necessary in the case the curl of $\mathbf{J}$ is equal to a non-zero Burgers vector. It is again obvious that for $a_c \to \infty$ we obtain $\varepsilon_{ijk} \partial_j J_k \to B_i$, i.e. the involution \eqref{burgers.involution}. 
However, the FO-CCZ4 system considered in the rest of this paper will only have homogeneous curl involutions, i.e. the additional field $\chi$ is not needed. }

{At this point we stress again very clearly that the proposed GLM curl cleaning approach comes at the \textit{expense} of needing to evolve one additional 3 vector $\psi_k$ and one additional scalar $\varphi$ for homogeneous curl involutions. For non-homogeneous curl involutions, one additional 3 vector $\psi_k$ and two additional cleaning scalars $\varphi$ and $\chi$ are needed. This is the price to pay for the great flexibility and ease of implementation of the proposed method. }

{We also would like to clearly stress at this point that the choice of the cleaning speeds is not unique and is left to the user. Compared to an exactly structure-preserving scheme, this may be considered as another major drawback of the proposed GLM cleaning approach, besides the large number of required additional evolution variables. }

{Note that in the special case of a \textit{linear} source term $S_k$ in \eqref{eqn.Jk.src}, e.g. 
\begin{equation} 
\partial_t J_k + \partial_k ( v_m J_m ) + v_m \left( \partial_m J_k - \partial_k J_m \right) =  -\frac{1}{\tau} J_k, \\ 
\label{eqn.Jk.linS} 
\end{equation} 	
with a \textit{constant} relaxation time $\tau>0$, it is obvious that the field $J_k$ will still remain curl-free and thus $B_k=0$ for
all times, if the curl of $J_k$ was initially zero. In other words, to generate a non-vanishing curl of $J_k$, one needs a \textit{nonlinear} source term, such as the one proposed in \cite{GPRmodel}, where the factor in front of $J_k$ in the source term also depends on temperature and density. 
}


\section{Governing equations of the FO-CCZ4 system with GLM curl cleaning} 
\label{sec.model} 
In compact 4D tensor notation, the Einstein field equations with matter source terms read 
\begin{equation} 
   R_{\mu \nu} - \frac{1}{2} g_{\mu \nu} R = 8 \pi T_{\mu \nu},
\end{equation} 
where $g_{\mu \nu}$ is the 4-metric of the spacetime, which is the primary unknown of the system, 
$R_{\mu \nu} = R^\lambda_{\, \, \mu \lambda \nu}$ is the 4-Ricci tensor, which is a 
contraction of the Riemann curvature tensor $R^\kappa_{\, \, \mu \lambda \nu}$ 
associated with the metric,  $R = g^{\mu \nu} R_{\mu \nu}$ is the 4-Ricci 
scalar and $T_{\mu \nu}$ is the energy-momentum tensor, which in this work will be considered
as an externally \textit{given} quantity. 
The original second-order CCZ4 governing system \cite{Alic:2012} is a nonlinear time-dependent 
PDE system with first order time derivatives and mixed first and second order spatial 
derivatives. It can be derived from the Z4 Lagrangian 
\begin{equation} 
\mathcal{L} = g^{\mu\nu} ( R_{\mu\nu} + 2 \nabla_\mu Z_\nu ),
\end{equation}
which adds a new vector $Z_\mu$ to the classical Einstein-Hilbert Lagrangian 
(see \cite{Bona2009} for details). The Z4 vector ($Z_\mu$ ) has been introduced in a series of papers by Bona et al. \cite{Bona:2002fq,Bona:2003fj,Bona:2003qn,Bona:2004yp} in order to enforce the nonlinear  
involutions of the Einstein field equations via a hyperbolic constraint-cleaning 
approach, which can be seen as an extension of the GLM method of Munz et al.  \cite{MunzCleaning,Dedneretal} to nonlinear involution constraints. In the Einstein field equations, 
the involutions are the so-called Hamiltonian constraint $\mathcal{H}$ and the momentum constraints $\mathcal{M}_i$ defined later. 
Additional algebraic constraint-damping terms can be added \cite{Gundlach2005:constraint-damping}, 
so that the Einstein field equations with Z4 constraint cleaning and constraint damping 
\cite{Alic:2009,Alic:2012} finally read 
\begin{equation}
R_{\mu\nu} + \nabla_{(\mu} Z_{\nu)} +
\kappa_1 \left( n_{(\mu} Z_{\nu)}
-(1+\kappa_2) g_{\mu\nu} n_\alpha Z^\alpha \right)
= 8 \pi \left( T_{\mu \nu} - \frac{1}{2} g_{\mu \nu} T \right), 
\end{equation}
where $T = g^{\mu \nu} T_{\mu\nu}$ denotes the trace of the energy-momentum tensor,   
$\boldsymbol{n}$ is the unit vector normal to the spatial hypersurfaces and the $\kappa_i$ 
are adjustable constants that are related to the damping of the Z4 vector. 
Since the covariant forms of the field equations above are not directly suitable 
for numerical treatment, a usual 3+1 decomposition of spacetime is adopted, see, \eg  \cite{Alcubierre:2008,Baumgarte2010,Bona2009}. Hence, the line element is written as  
\begin{equation}
d s^2 = -\alpha^2 d t^2 
+ \gamma_{ij} \,  \left( d x^i + \beta^i d t \right) \left( d x^j + \beta^j d t \right),
\end{equation}
with the usual lapse function $\alpha$, the shift vector $\beta^i$ and the spatial 3-metric 
$\gamma_{ij}$, which is related to the 4-metric by $\gamma_{\mu \nu} = g_{\mu \nu} + n_\mu n_\nu$. 
The 3+1 split then leads to evolution equations for $\gamma_{ij}$, as well 
as for the extrinsic curvature $K_{ij} = -\frac{1}{2}\mathcal{L}_n \gamma_{ij}$, 
where $\mathcal{L}_n$ denotes the Lie derivative along the vector $n^\mu$. Due to the  
freedom concerning the choice of the coordinate system in general relativity (gauge freedom), 
the lapse and the shift can in principle be freely chosen. The four constraint equations 
of the ADM system, namely the Hamiltonian constraint 
\begin{equation}
\mathcal{H} = R_{ij} g^{ij} - K_{ij} K^{ij} + K^2 - 16 \pi \tau = 0
\label{eqn.ham} 
\end{equation}
and the three momentum constraints 
\begin{equation}
\mathcal{M}_i = \gamma^{jl} \left( \partial_l K_{ij} - \partial_i K_{jl} - \Gamma^m_{\,j,l} K_{mi} + \Gamma^m_{\,ji} K_{ml} \right) - 8 \pi S_i  = 0,  
\label{eqn.mom} 
\end{equation}
including the matter source terms, are taken into account by additional evolution equations for the $Z_\mu$ vector. In \eqref{eqn.ham} and \eqref{eqn.mom} $\tau = n^\mu n^\nu T_{\mu \nu}$ and $S_i = -P_{i}^{\ \mu} n^\nu T_{\mu \nu}$ (with the spatial projector $P_\nu^{\ \mu} = \delta_{\nu}^{\ \mu} + n_\nu n^{\mu}$) are projections of the energy-momentum tensor into the spatial hyperplane and represent 
the usual evolution quantities for total energy and momentum in the general relativistic Euler equations \cite{DelZanna2002}. 
The CCZ4 formulation of \cite{Alic:2012} introduces a conformal factor $\phi := |\gamma_{ij}|^{-1/6}$ in order to define a conformal 3-metric with unit determinant as 
\begin{equation}
 \tilde\gamma_{ij}:= \phi^2\gamma_{ij}, \qquad |\tilde\gamma_{ij}|=1. 
\end{equation}
Like in the BSSNOK system \cite{Shibata95,Baumgarte99,Nakamura87}, the extrinsic curvature is   
decomposed into its trace-free part $\tilde{A}_{ij}$ 
\begin{equation}
\tilde A_{ij} := \phi^2\left(K_{ij}-\frac{1}{3}K\gamma_{ij}\right), 
\end{equation}
and into its trace $K=K_{ij}\gamma^{ij}$, which both become new evolution variables. More details about the original second order CCZ4 system and its derivation can be found in \cite{Alic:2012},
while a detailed discussion of its first order reduction together with a proof of strong hyperbolicity were given in \cite{ADERCCZ4}.  

In the following, we show how our new GLM curl cleaning approach described on the toy system in the previous section can be applied to the full FO-CCZ4 formulation of the Einstein field equations. 
The FO-CCZ4 system as proposed in \cite{ADERCCZ4} has a very particular split structure, where the  quantities that define the 4-metric, namely the lapse function $\alpha$, the shift vector $\beta^i$, the conformal spatial metric tensor $\tilde{\gamma}_{ij}$ and the conformal factor $\phi$ are  governed by an ODE-type system, i.e. by a system that contains only first order time derivatives and purely algebraic source terms, but no spatial derivatives:  
\begin{subequations}
	\begin{eqnarray}
	\label{eqn.gamma}
	\partial_t\tilde\gamma_{ij}  
	&=&  {\beta^k 2 D_{kij} + \tilde\gamma_{ki} B_{j}^k  + \tilde\gamma_{kj} B_{i}^k - \frac{2}{3}\tilde\gamma_{ij} B_k^k }
	- 2\alpha \left( \tilde A_{ij} - {\frac{1}{3} \tilde \gamma_{ij} {\rm tr}{\tilde A} } \right),
	\\
	\label{eqn.alpha}
	{ \partial_t \ln{\alpha} }  &=&  { \beta^k A_k } - \alpha g(\alpha) ( K - K_0 - 2\Theta {c} ), 
	\qquad \qquad \partial_t K_0 = 0, \\
	\label{eqn.beta}
	\partial_t \beta^i   &=& 
	s \beta^k B_k^i +
	s f b^i, \\
	\label{eqn.phi}
	{ \partial_t \ln{\phi} }  &=&  { \beta^k P_k } + \frac{1}{3} \left( \alpha K - {B_k^k} \right),
	\end{eqnarray} 
\end{subequations}
The use of the logarithms in the evolution equations for the lapse and the conformal factor is for convenience, in order to always guarantee \textit{positivty} of $\alpha$ and $\phi$ also at the 
\textit{discrete level}.  
The next primary evolution quantities are the trace-free extrinsic curvature tensor $\tilde A _{ij}$,
the trace of the extrinsic curvature $K$, the GLM cleaning scalar $\Theta$ that accounts for the Hamiltonian constraint $\mathcal{H}$, the modified contracted Christoffel symbols $\hat\Gamma^i$, which contain the spatial part of the Z4 cleaning vector $Z^i = \frac{1}{2} \phi^2 \left( \hat\Gamma^i - \tilde \Gamma^i \right)$ for the cleaning of the momentum constraints, and the variable $b^i$ that is needed for the gamma driver shift condition:   
\begin{subequations}
	\begin{eqnarray}
	\label{eqn.Aij}
	\partial_t\tilde A _{ij} - \beta^k \partial_k\tilde A_{ij}  + \phi^2 \bigg[ \nabla_i\nabla_j  \alpha - \alpha \left( R_{ij}+ \nabla_i Z_j + \nabla_j Z_i - 8 \pi S_{ij} \right) \bigg] 	 
	-  \frac{1}{3} \tilde\gamma_{ij} \bigg[ \nabla^k \nabla_k \alpha - \alpha (R +  2 \nabla_k Z^k - 8 \pi S ) \bigg] &&    \nonumber \\	
	 =  { \tilde A_{ki} B_j^k + \tilde A_{kj} B_i^k - \frac{2}{3}\tilde A_{ij} B_k^k }
	+ \alpha \tilde A_{ij}(K - 2 \Theta {c} ) - 2 \alpha\tilde A_{il} \tilde\gamma^{lm} \tilde A_{mj}, \ && 
	\\ \nonumber \\  
	\label{eqn.K}
	\partial_t K - \beta^k \partial_k K + \nabla^i \nabla_i \alpha - \alpha( R + 2 \nabla_i Z^i) =
	\alpha K (K - 2\Theta {c} ) - 3\alpha\kappa_1(1+\kappa_2)\Theta + 4 \pi \alpha( S - 3 \tau), \ && 
	\\ \nonumber \\  
	\label{eqn.theta}
	\partial_t \Theta - \beta^k\partial_k\Theta - \frac{1}{2}\alpha {e^2} ( R + 2 \nabla_i Z^i)
	=  \alpha {e^2} \left( \frac{1}{3} K^2 - \frac{1}{2} \tilde{A}_{ij} \tilde{A}^{ij}  - 8 \pi \tau \right) - \alpha \Theta K {c} - {Z^i \alpha A_i}
	- \alpha\kappa_1(2+ \kappa_2)\Theta, \  &&   
	\\ \nonumber \\    
	\label{eqn.Ghat}
	\partial_t \hat\Gamma^i - \beta^k \partial_k \hat \Gamma^i + \frac{4}{3} \alpha \tilde{\gamma}^{ij} \partial_j K  - 2 \alpha \tilde{\gamma}^{ki} \partial_k \Theta
	- \tilde{\gamma}^{kl} \partial_{(k} B_{l)}^i
	- \frac{1}{3} \tilde{\gamma}^{ik}  \partial_{(k} B_{l)}^l - { s 2 \alpha \tilde{\gamma}^{ik}  \tilde{\gamma}^{nm} \partial_k \tilde{A}_{nm}   }
	&& \nonumber \\
     =  { \frac{2}{3} \tilde{\Gamma}^i B_k^k - \tilde{\Gamma}^k B_k^i  } +
	2 \alpha \left( \tilde{\Gamma}^i_{jk} \tilde{A}^{jk} - 3 \tilde{A}^{ij} P_j \right) -
	2 \alpha \tilde{\gamma}^{ki} \left( \Theta A_k + \frac{2}{3} K Z_k \right) -
	2 \alpha \tilde{A}^{ij} A_j && \nonumber \\
	  - { 4 s \, \alpha \tilde{\gamma}^{ik} D_k^{~\,nm} \tilde{A}_{nm} } + 2\kappa_3 \left( \frac{2}{3} \tilde{\gamma}^{ij} Z_j B_k^k - \tilde{\gamma}^{jk} Z_j B_k^i \right) - 2 \alpha \kappa_1 \tilde{\gamma}^{ij} Z_j - 16 \pi \alpha \tilde{\gamma}^{ij} S_j,  \   &&
	  \\  
	\label{eqn.b}
	\partial_t b^i - s \beta^k \partial_k b^i = s \left(  \partial_t \hat\Gamma^i - \beta^k \partial_k \hat \Gamma^i - \eta b^i \right), && 
	\end{eqnarray} 
\end{subequations}
In order to obtain a strongly hyperbolic first order reduction, the following evolution system 
for the auxiliary variables 
$A_k := \partial_k \alpha / \alpha $, $B_k^{i} := \partial_k\beta^i$, 
$D_{kij} := \frac{1}{2}\partial_k\tilde\gamma_{ij}$ and 
$P_k :=  \partial_k \phi / \phi$ is added: 
\begin{subequations}
	\begin{eqnarray}
		\label{eqn.A}
	\partial_t A_{k} - {\beta^l \partial_l A_k} + \alpha g(\alpha) \left( \partial_k K - \partial _k K_0 - 2c \partial_k \Theta \right)
	+ {s \alpha g(\alpha) \tilde{\gamma}^{nm} \partial_k \tilde{A}_{nm} }\\ \nonumber  =  
	+ {2s\, \alpha g(\alpha) D_k^{~\,nm} \tilde{A}_{nm} }
	-\alpha A_k \left( K - K_0 - 2 \Theta c \right) h(\alpha) + B_k^l ~A_{l} \,,
	\\ \nonumber \\ 
	\label{eqn.B}
	\partial_t B_k^i  - s\beta^l \partial_l B_k^i - s\left(  f \partial_k b^i - { \alpha^2 \mu \, \tilde{\gamma}^{ij} \left( \partial_k P_j - \partial_j P_k \right)
		+ \alpha^2 \mu \, \tilde{\gamma}^{ij} \tilde{\gamma}^{nl} \left( \partial_k D_{ljn} - \partial_l D_{kjn} \right) } \right)
	= 
	s B^l_k~B^i_l \,,
	\\ \nonumber \\ 
	\label{eqn.D}
	\partial_t D_{kij}  - {\beta^l \partial_l D_{kij}} + s \left(
	- \frac{1}{2} \tilde{\gamma}_{mi} \partial_{(k} {B}_{j)}^m
	- \frac{1}{2} \tilde{\gamma}_{mj} \partial_{(k} {B}_{i)}^m
	+ \frac{1}{3} \tilde{\gamma}_{ij} \partial_{(k} {B}_{m)}^m  \right)
	+  \alpha \partial_k \tilde{A}_{ij}
	-  { \alpha \frac{1}{3} \tilde{\gamma}_{ij} \tilde{\gamma}^{nm} \partial_k \tilde{A}_{nm} }   
	\nonumber \\ 
	=   B_k^l D_{lij} + B_j^l D_{kli} + B_i^l D_{klj} - \frac{2}{3} B_l^l D_{kij} - {  \alpha \frac{2}{3}  \tilde{\gamma}_{ij} D_k^{~\,nm} \tilde{A}_{nm} } - \alpha A_k \left( \tilde{A}_{ij} - \frac{1}{3} \tilde{\gamma}_{ij} {\rm tr} \tilde{A} \right), 
	 \\ \nonumber \\ 
	\label{eqn.P}
	\partial_t P_{k} - \beta^l \partial_l P_{k} - \frac{1}{3} \alpha \partial_k K
	+ \frac{1}{3} \partial_{(k} {B}_{i)}^i   - {s \frac{1}{3} \alpha \tilde{\gamma}^{nm} \partial_k \tilde{A}_{nm} }  = 
	\frac{1}{3} \alpha A_k K + B_k^l P_l - {s \frac{2}{3} \alpha \, D_k^{~\,nm} \tilde{A}_{nm} }. 
	\end{eqnarray}
\end{subequations}
The governing PDE system \eqref{eqn.gamma}-\eqref{eqn.P} contains the following terms as defined below:  
\begin{eqnarray}
{\rm tr} \tilde{A} & = &  \tilde{\gamma}^{ij} \tilde{A}_{ij}, \qquad \textnormal{ and } \qquad \tilde{\gamma} = \textnormal{det}( \tilde{\gamma}_{ij} ), \\
\partial_k \tilde{\gamma}^{ij} & = &  - 2 \tilde{\gamma}^{in} \tilde{\gamma}^{mj} D_{knm} := -2 D_k^{~\,ij}, 
\qquad\text{(theorem of the derivative of the inverse matrix)} \\
\tilde{\Gamma}_{ij}^k &=& \tilde{\gamma}^{kl} \left( D_{ijl} + D_{jil} - D_{lij} \right), \\
\label{eqn.dchr}
\partial_k \tilde{\Gamma}_{ij}^m & = & -2 D_k^{ml} \left( D_{ijl} + D_{jil} - D_{lij} \right)
+ \tilde{\gamma}^{ml} \left( \partial_{(k} {D}_{i)jl} + \partial_{(k} {D}_{j)il} - \partial_{(k} {D}_{l)ij} \right) , \\ %
\Gamma_{ij}^k &=& \tilde{\gamma}^{kl} \left( D_{ijl} + D_{jil} - D_{lij} \right) - \tilde{\gamma}^{kl} \left( \tilde{\gamma}_{jl} P_i + \tilde{\gamma}_{il} P_j - \tilde{\gamma}_{ij} P_l \right)
= \tilde{\Gamma}_{ij}^k - \tilde{\gamma}^{kl} \left( \tilde{\gamma}_{jl} P_i + \tilde{\gamma}_{il} P_j - \tilde{\gamma}_{ij} P_l \right),  \\
\partial_k \Gamma_{ij}^m &=& -2 D_k^{ml} \left( D_{ijl} + D_{jil} - D_{lij} \right) +
2 D_k^{ml} \left( \tilde{\gamma}_{jl} P_i + \tilde{\gamma}_{il} P_j  - \tilde{\gamma}_{ij} P_l \right)
- 2 \tilde{\gamma}^{ml} \left(  D_{kjl} P_i + D_{kil} P_j  - D_{kij} P_l \right)
\nonumber \\
&&	   + \tilde{\gamma}^{ml} \left( \partial_{(k} {D}_{i)jl} + \partial_{(k} {D}_{j)il} - \partial_{(k} {D}_{l)ij} \right)  -
\tilde{\gamma}^{ml} \left( \tilde{\gamma}_{jl} \partial_{(k} {P}_{i)} + \tilde{\gamma}_{il} \partial_{(k} {P}_{j)}  - \tilde{\gamma}_{ij} \partial_{(k} {P}_{l)} \right) , \\
R^m_{ikj} & = & \partial_k \Gamma^m_{ij} - \partial_j \Gamma^m_{ik} + \Gamma^l_{ij} \Gamma^m_{lk} - \Gamma^l_{ik} \Gamma^m_{lj}, \\
R_{ij} & = &  R^m_{imj}, \\
\nabla_i \nabla_j \alpha &=& \alpha A_i A_j - \alpha \Gamma^k_{ij} A_k + \alpha \partial_{(i} {A}_{j)}, \\
\nabla^i \nabla_i \alpha &=& \phi^2 \tilde{\gamma}^{ij} \left( \nabla_i \nabla_j \alpha \right), \\
\tilde{\Gamma}^i & = &  \tilde{\gamma}^{jl} \tilde{\Gamma}^i_{jl},  \\
\partial_k \tilde{\Gamma}^i & = & -2 D_k^{jl} \, \tilde{\Gamma}^i_{jl} + \tilde{\gamma}^{jl} \, \partial_k \tilde{\Gamma}^i_{jl}, \\
Z_i &=& \frac{1}{2} \tilde\gamma_{ij} \left(\hat{\Gamma}^j-\tilde{\Gamma}^j \right), \qquad  Z^i = \frac{1}{2} \phi^2 (\hat\Gamma^i-\tilde\Gamma^i), \\
\nabla_i Z_j &=& D_{ijl} \left(\hat{\Gamma}^l-\tilde{\Gamma}^l \right) + \frac{1}{2} \tilde\gamma_{jl} \left( \partial_i \hat{\Gamma}^l - \partial_i \tilde{\Gamma}^l \right) - \Gamma^l_{ij} Z_l, \\
R + 2 \nabla_k Z^k & = & \phi^2 \tilde{\gamma}^{ij} \left( R_{ij} +
\nabla_i Z_j + \nabla_j Z_i \right)\,, \\
h(\alpha) &=& \left( g(\alpha) + \alpha \frac{\partial g(\alpha)}{ \partial \alpha}  \right).
%
%
%
%
\end{eqnarray}
The function $g(\alpha)$ in the PDE for the lapse $\alpha$ controls the slicing condition, where
$g(\alpha)=1$ leads to harmonic slicing and $g(\alpha)=2/\alpha$ leads to the so-called $1+\log$ 
slicing condition, see \cite{Bona95b}.  
As already mentioned above, the auxiliary quantities $A_k$, $P_k$, $B^i_k$ and $D_{kij}$ are defined as (scaled) spatial gradients of the primary variables 
$\alpha$, $\phi$, $\beta^i$ and $\tilde{\gamma}_{ij}$, respectively, and read: 
\begin{equation}
\label{eq:Auxiliary}
A_i := \partial_i\ln\alpha = \frac{\partial_i \alpha }{\alpha}\,, \qquad
B_k^{i} := \partial_k\beta^i\,,
\qquad
D_{kij} := \frac{1}{2}\partial_k\tilde\gamma_{ij}\,, \qquad
P_i       := \partial_i\ln\phi = \frac{\partial_i \phi}{\phi}\,.
\end{equation}
Hence, as a result, they must satisfy the following curl involutions or so-called second order ordering constraints \cite{Alic:2009,Gundlach:2005ta}:  
\begin{align}
\label{eqn.second.ord.const}
\mathcal{A}_{lk}   &:= \partial_l A_k     - \partial_k A_l     = 0, & 
\mathcal{P}_{lk}   &:= \partial_l P_k     - \partial_k P_l     = 0, \nonumber \\
\mathcal{B}_{lk}^i &:= \partial_l B_k^i   - \partial_k B_l^i   = 0, & 
\mathcal{D}_{lkij} &:= \partial_l D_{kij} - \partial_k D_{lij} = 0.   
\end{align}
In the governing PDE system above, we have already made use of these curl involutions by \textit{symmetrizing} 
the spatial derivatives of the auxiliary variables as follows:
\begin{equation}
\partial_{(k} {A}_{i)}     := \frac{ \partial_k A_i + \partial_i A_k       }{2}, \quad
\partial_{(k} {P}_{i)}     := \frac{ \partial_k P_i + \partial_i P_k       }{2}, \quad
\partial_{(k} {B}^i_{j)}   := \frac{ \partial_k B^i_j + \partial_j B^i_k     }{2}, \quad
\partial_{(k} {D}_{l)ij}  := \frac{ \partial_k D_{lij} + \partial_l D_{kij} }{2}.
\label{eqn.symm.aux}
\end{equation} 
The new hyperbolic curl cleaning approach applied to system \eqref{eqn.gamma}-\eqref{eqn.P} obviously affects only the governing PDE for the auxiliary variables $A_k$, $P_k$, $B^i_k$ and $D_{kij}$, which are now \textit{augmented} as follows: 
\begin{subequations}
\begin{eqnarray}
\label{eqn.A.GLM}
\partial_t A_{k} - {\beta^l \partial_l A_k} \textcolor{blue}{ \, + \, \varepsilon_{klm} \, \partial_l \psi^A_m } &+& \alpha g(\alpha) \left( \partial_k K - \partial _k K_0 - 2c \partial_k \Theta \right)
+ {s \alpha g(\alpha) \tilde{\gamma}^{nm} \partial_k \tilde{A}_{nm} } 
 \nonumber \\   
 &=&  
+ {2s\, \alpha g(\alpha) D_k^{~\,nm} \tilde{A}_{nm} }
-\alpha A_k \left( K - K_0 - 2 \Theta c \right) h(\alpha) + B_k^l ~A_{l} \,, \\  
\textcolor{blue}{\partial_t \psi^A_k - (a_c^A)^2 \, \varepsilon_{klm} \, \partial_l A_m }  
\textcolor{red}{ \, + \, \partial_k \varphi^A } & = & \textcolor{blue}{-\epsilon^A_c \psi^A_k}, \\ 
\textcolor{red}{\partial_t \varphi^A + (a_d^A)^2 \, \partial_m \psi^A_m } 
 & = & \textcolor{red}{-\epsilon^A_d \varphi^A}, 
 \\ \nonumber \\  
\label{eqn.B.GLM}
\partial_t B_k^i - s\beta^l \partial_l B_k^i \textcolor{blue}{ + \varepsilon_{klm} \, \partial_l (\psi^B)^i_m } &-& s\left(  f \partial_k b^i - { \alpha^2 \mu \, \tilde{\gamma}^{ij} \left( \partial_k P_j - \partial_j P_k \right)
	+ \alpha^2 \mu \, \tilde{\gamma}^{ij} \tilde{\gamma}^{nl} \left( \partial_k D_{ljn} - \partial_l D_{kjn} \right) } \right) \nonumber \\ 
   &=& s B^l_k~B^i_l, \\ 
\textcolor{blue}{\partial_t (\psi^B)^i_k - (a_c^B)^2 \, \varepsilon_{klm} \, \partial_l B^i_m } 
\textcolor{red}{\, + \, \partial_k (\varphi^B)^i } & = & \textcolor{blue}{-\epsilon^B_c (\psi^B)^i_k}, \\ 
\textcolor{red}{\partial_t (\varphi^B)^i + (a_d^B)^2 \, \partial_m (\psi^B)^i_m } 
& = & \textcolor{red}{-\epsilon^B_d (\varphi^B)^i },  
\\ \nonumber \\  
\label{eqn.D.GLM}
\partial_t D_{kij} - {\beta^l \partial_l D_{kij}} \textcolor{blue}{ \, + \, \varepsilon_{klm} \, \partial_l \psi^D_{mij} } &+& s \left(
- \frac{1}{2} \tilde{\gamma}_{mi} \partial_{(k} {B}_{j)}^m
- \frac{1}{2} \tilde{\gamma}_{mj} \partial_{(k} {B}_{i)}^m
+ \frac{1}{3} \tilde{\gamma}_{ij} \partial_{(k} {B}_{m)}^m  \right)  \nonumber \\ 
&+&  \alpha \partial_k \tilde{A}_{ij} 
 -   { \alpha \frac{1}{3} \tilde{\gamma}_{ij} \tilde{\gamma}^{nm} \partial_k \tilde{A}_{nm} }   
\nonumber  \\ 
&=&   B_k^l D_{lij} + B_j^l D_{kli} + B_i^l D_{klj} - \frac{2}{3} B_l^l D_{kij} - {  \alpha \frac{2}{3}  \tilde{\gamma}_{ij} D_k^{~\,nm} \tilde{A}_{nm} } \nonumber \\ 
&& - \alpha A_k \left( \tilde{A}_{ij} - \frac{1}{3} \tilde{\gamma}_{ij} {\rm tr} \tilde{A} \right), \\
\textcolor{blue}{\partial_t \psi^D_{kij} - (a_c^D)^2 \, \varepsilon_{klm} \, \partial_l D_{mij} } 
\textcolor{red}{ \, + \, \partial_k \varphi^D_{ij} } & = & \textcolor{blue}{-\epsilon^D_c \psi^D_{kij} } \\ 
\textcolor{red}{\partial_t \varphi^D_{ij} + (a_d^D)^2 \, \partial_m \psi^D_{mij} } 
& = & \textcolor{red}{-\epsilon^D_d \varphi^D_{ij} }   
\\ \nonumber \\  
\label{eqn.P.GLM}
\partial_t P_{k} - \beta^l \partial_l P_{k} \textcolor{blue}{ \, + \, \varepsilon_{klm} \, \partial_l \psi^P_{m} } &-& \frac{1}{3} \alpha \partial_k K
+ \frac{1}{3} \partial_{(k} {B}_{i)}^i   - {s \frac{1}{3} \alpha \tilde{\gamma}^{nm} \partial_k \tilde{A}_{nm} } \\  \nonumber &=&
\frac{1}{3} \alpha A_k K + B_k^l P_l - {s \frac{2}{3} \alpha \, D_k^{~\,nm} \tilde{A}_{nm} }. 
\\  
\textcolor{blue}{\partial_t \psi^P_k - (a_c^P)^2 \, \varepsilon_{klm} \, \partial_l P_m }  
\textcolor{red}{ \, + \, \partial_k \varphi^P } & = & \textcolor{blue}{-\epsilon^P_c \psi^P_k} \\ 
\textcolor{red}{\partial_t \varphi^P + (a_d^P)^2 \, \partial_m \psi^P_m } 
& = & \textcolor{red}{-\epsilon^P_d \varphi^P}. 
\label{eqn.GLM.end} 
\end{eqnarray}
\end{subequations}
Here, we have used essentially the same notation for the cleaning fields as the one employed in the previous section for the toy system. We have again used blue color in order to highlight the additional terms that are responsible for the curl cleaning and in red those that are used for the divergence cleaning. The associated cleaning fields are $\psi^A_k$, $\psi^P_k$, $(\psi^B)_k^i$ and 
$\psi^D_{kij}$ and $\varphi^A$, $\varphi^P$, $(\varphi^B)_k^i$ and $\varphi^D_{kij}$ for the
variables $A_k$, $P_k$, $B_k^i$ and $D_{kij}$, respectively. 
The final FO-CCZ4 system with GLM cleaning is therefore given by \eqref{eqn.gamma}-\eqref{eqn.phi}, \eqref{eqn.Aij}-\eqref{eqn.b} and \eqref{eqn.A.GLM}-\eqref{eqn.GLM.end}, i.e. the set of equations
\eqref{eqn.A.GLM}-\eqref{eqn.GLM.end} replaces the equations for the auxiliary variables \eqref{eqn.A}-\eqref{eqn.P} of the original FO-CCZ4 system. 

{At this point we clearly emphasize that the additional equations for the GLM cleaning inside the proposed augmented FO-CCZ4 system are not covariant. However, in the view of the authors this is not a problem, since all these additional equations have no physical  meaning but are only needed to properly account for the involution constraints on the discrete level when applying a numerical scheme to \eqref{eqn.gamma}-\eqref{eqn.GLM.end}. For the same reason, it is also possible to choose superluminal cleaning speeds, in particular because the involutions are asymptotically retrieved in the limit when the cleaning speeds tend to infinity. Of course, the negative side effect of superluminal cleaning speeds is an increased numerical dissipation for the physical quantities and an increased computational effort, due to the CFL stability condition of explicit schemes. } 

%
\section{Numerical results}       
\label{sec.results}
The augmented FO-CCZ4 system with GLM curl cleaning presented in the previous section can formally 
be written as one big first order hyperbolic PDE system of 103 evolution quantities: 
\begin{equation}
   \partial_t \Q + \A(\Q) \cdot \nabla \Q= \S(\Q), 
   \label{eqn.pde} 
\end{equation}
where $\Q = (\alpha, \beta^i, \tilde \gamma_{ij}, \phi, K_0, \tilde A_{ij}, K, \Theta, \hat \Gamma^i, b^i, A_k, \psi^A_k, \varphi^A, B_k^i, (\psi^B)^i_k, (\varphi^B)^i, D_{kij}, \psi^D_{kij}, \varphi^D_{ij}, P_k,
\psi^P_k, \varphi^P )^T$
is the state vector, $\A(\Q) = (\A_1, \A_2, \A_3)$ are the system matrices in the three coordinate
directions and $\S(\Q)$ is the algebraic source term on the right hand side. For the complete 
set of eigenvalues and associated eigenvectors of the original FO-CCZ4 system without GLM cleaning, 
see \cite{ADERCCZ4}. The governing PDE system \eqref{eqn.pde} is now solved numerically with the aid of 
a high order accurate fully-discrete one-step ADER discontinuous Galerkin (DG) finite element scheme, 
exactly as described in \cite{DGLimiter1,Zanotti2015b,Zanotti2015c,GPRmodel,GPRmodelMHD,ADERCCZ4,Axioms}. 
Since all the technical details of the ADER-DG scheme can be found in these references, we omit the
description of the numerical method here, since this is not the main focus of this paper. Indeed, \textit{any} suitable discretization for hyperbolic systems of the type \eqref{eqn.pde} could be used in principle, and this \textit{total independence} of the underlying numerical method and grid topology is indeed the main strength of the GLM approach. In the following, the performance of the new GLM curl cleaning approach proposed in this paper is assessed on several benchmark problems. In all the following tests we use uniform 
Cartesian grids and set $K_0=0$. 

\subsection{Robust stability test} 
\label{sec.robstab} 

The first test case under consideration is the so-called robust stability test, which is a standard 
test problem in numerical general relativity and is taken from \cite{Alcubierre:2003pc}. 
The three-dimensional computational domain is given by the box $\Omega = [-0.5,0.5]^3$.
In this test problem, matter is absent, hence $\tau = 0$ and $S_i = 0$.  

As gauge conditions we employ a frozen shift condition $\partial_t
\beta^i = 0$ by setting $s=0$ in the FO-CCZ4 system, together with a harmonic lapse,
which corresponds to $g(\alpha)=1$. The GLM cleaning speeds are set to $e=2$ for the 
cleaning of the nonlinear ADM constraints, while the cleaning speeds for the curl involutions
are chosen as $a_c^A = a_c^P = a_c^D = 1.5$, $a_d^A = a_d^D = a_d^P = 2$, and the respective 
damping parameters are set to $\epsilon_{c,d}^A =  \epsilon_{c,d}^P = \epsilon_{c,d}^D= 1$. 
The remaining constants in the FO-CCZ4 system are chosen as $\kappa_1=\kappa_2=\kappa_3=0$,
$c=0$ and $\eta=0$. As usual, in this test we start from a flat Minkowski metric 
\begin{equation}
d s^2 = -d t^2 + d x^2 + d y^2 + d z^2.
\label{eqn.robstab.metric}
\end{equation}
Different simulations are run with an unlimited ADER-DG $P_3$ scheme (polynomial approximation degree $N=3$)
on two different meshes with refinement factor $\rho \in \left\{ 1, 2 \right\}$. 
The meshes are composed of $(10 \cdot \rho)^3$ elements, corresponding to  
$((N+1) \cdot 10 \cdot \rho)^3$ spatial degrees of freedom. 
We then add uniformly distributed \textit{random perturbations} with perturbation 
amplitude $\epsilon = 10^{-6}/\rho^2$ to \textit{all} variables of the FO-CCZ4 system.   
Note that the chosen perturbations are four orders of magnitude larger that those suggested 
in \cite{Alcubierre:2003pc}. 

\begin{figure}
	\centering
	\begin{tabular}{cc} 
	\includegraphics[width=0.45\textwidth]{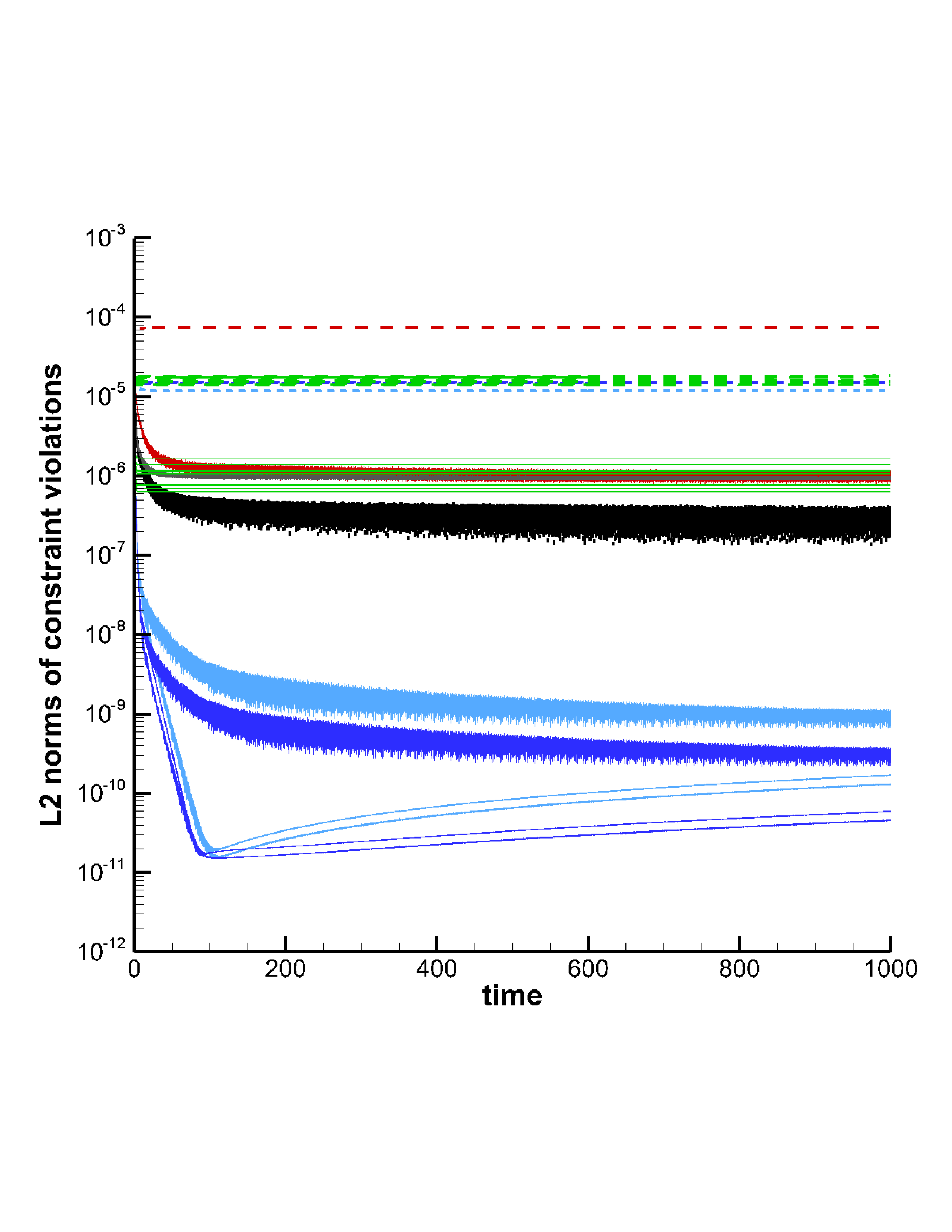}   & 
	\includegraphics[width=0.45\textwidth]{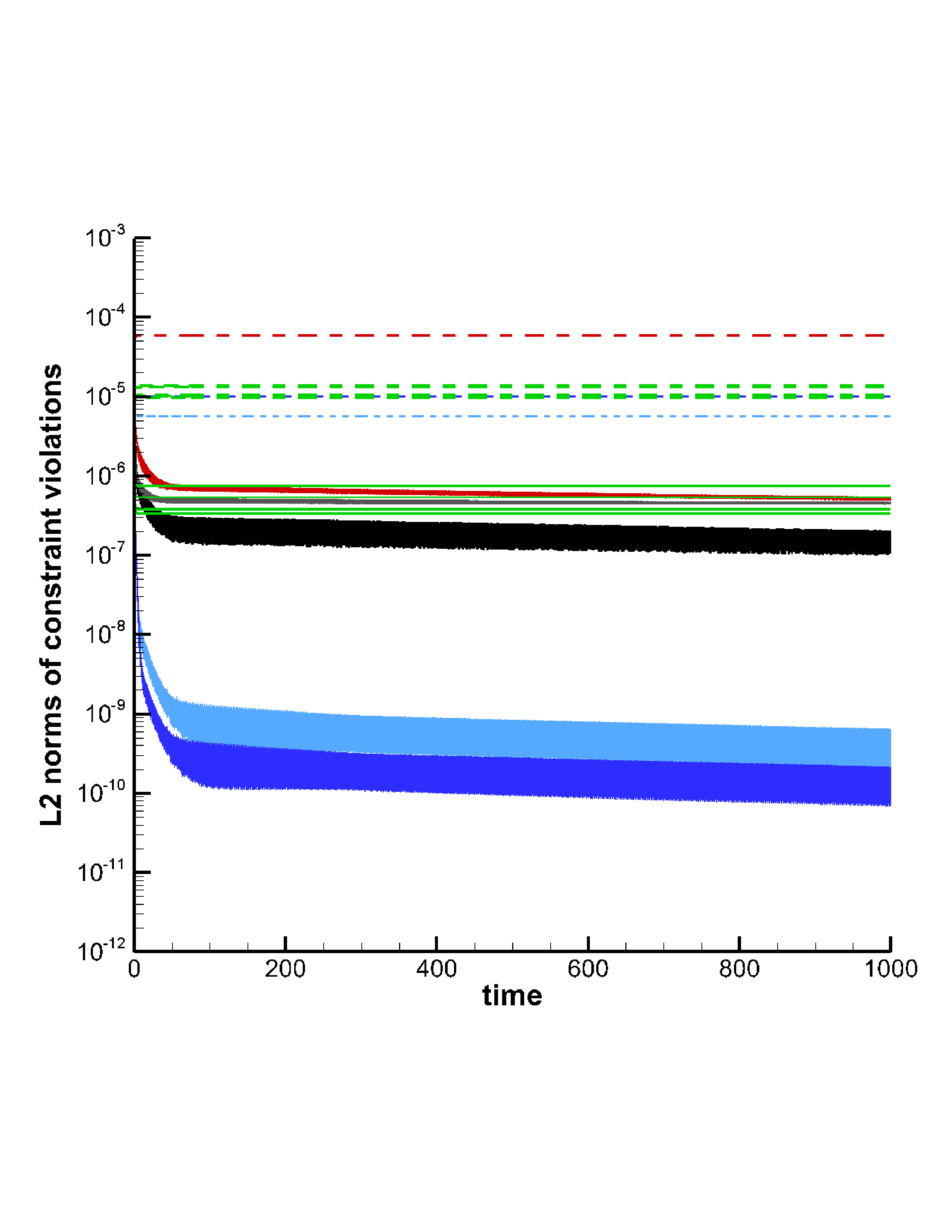}     
	\end{tabular} 
	
	\caption{Long-time constraint evolution for the robust stability test on the coarse mesh (left) and on the fine mesh (right). Dashed lines: FO-CCZ4 system without curl cleaning. Continuous lines: FO-CCZ4 with the new hyperbolic GLM curl cleaning. 
			 The following constraints are monitored: Hamiltonian constraint $\mathcal{H}$ (red), momentum constraints $\mathcal{M}_i$ (black), curl constraints $\mathcal{A}_{lk}$ (light blue), $\mathcal{P}_{lk}$ (dark blue) and $\mathcal{D}_{lkij}$ (green). }
	\label{fig.robstab}
\end{figure}

The time evolution of the four ADM constraints and of all curl involutions until $t=1000$ is 
reported in Fig. \ref{fig.robstab} for both simulations. For comparison, we also show the 
results obtained with the original FO-CCZ4 system \cite{ADERCCZ4} without curl cleaning. 
One can observe that the new hyperbolic GLM curl cleaning proposed in this paper reduces the 
errors in the Hamiltonian constraint by two orders of magnitude. The curl constraints for the 
auxiliary variable $D_{kij}$ improve by one order of magnitude, while the involutions 
$\mathcal{A}_{lk}$ and $\mathcal{P}_{lk}$ improve by more than \textit{four orders of magnitude}. 
These results clearly show the effectiveness of the new GLM curl cleaning approach proposed
in this paper for the FO-CCZ4 formulation of the Einstein field equations.  

\subsection{Stable neutron star} 
\label{sec.tovstar} 

\begin{figure}
	\centering
		\includegraphics[width=0.75\textwidth]{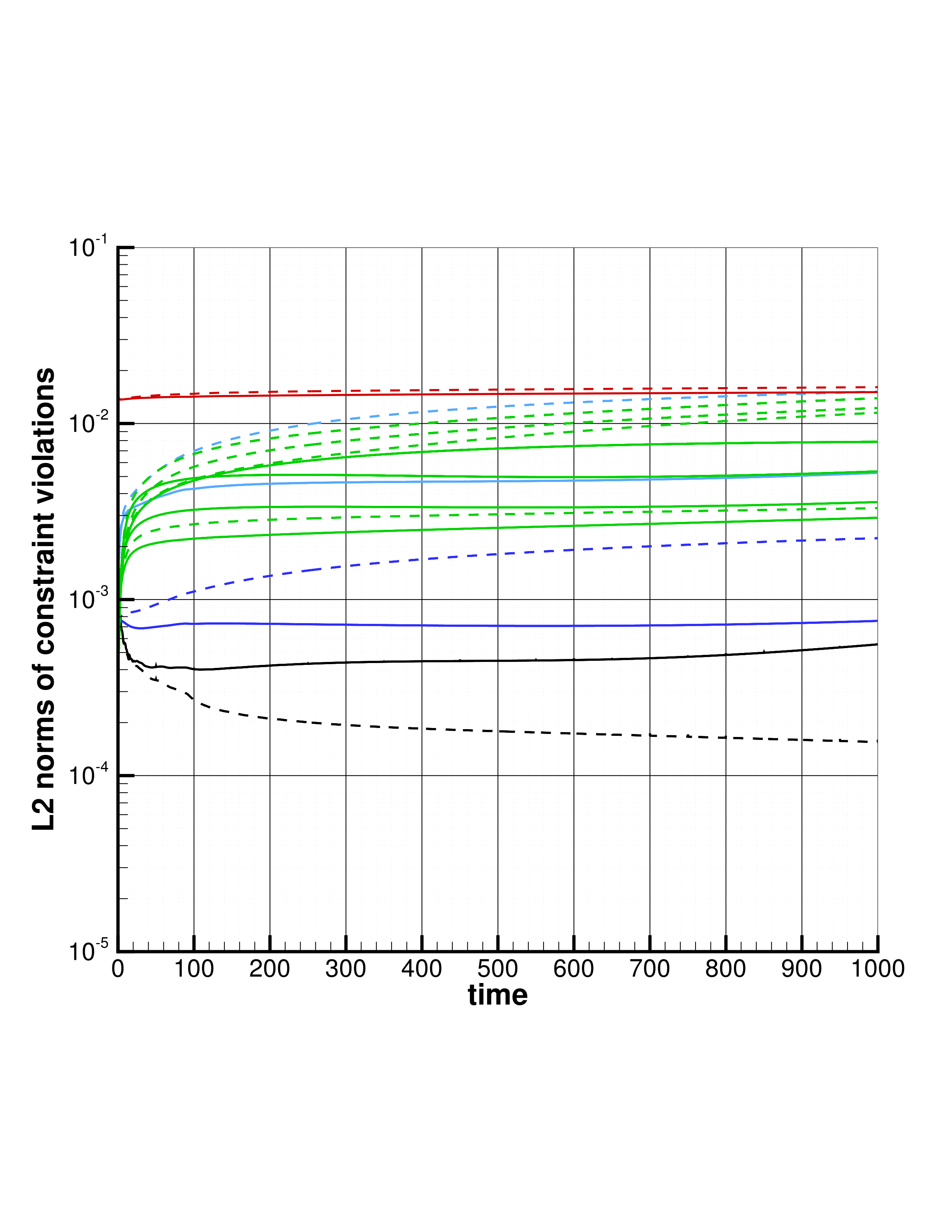}   
	\caption{Long-time constraint evolution for a stable neutron star until $t=1000$. Dashed lines: FO-CCZ4 system without curl cleaning. Continuous lines: FO-CCZ4 with hyperbolic GLM curl cleaning. }
	\label{fig.constraints}
\end{figure}

\begin{figure}
	\centering
	\begin{tabular}{cc}
		\includegraphics[width=0.45\textwidth]{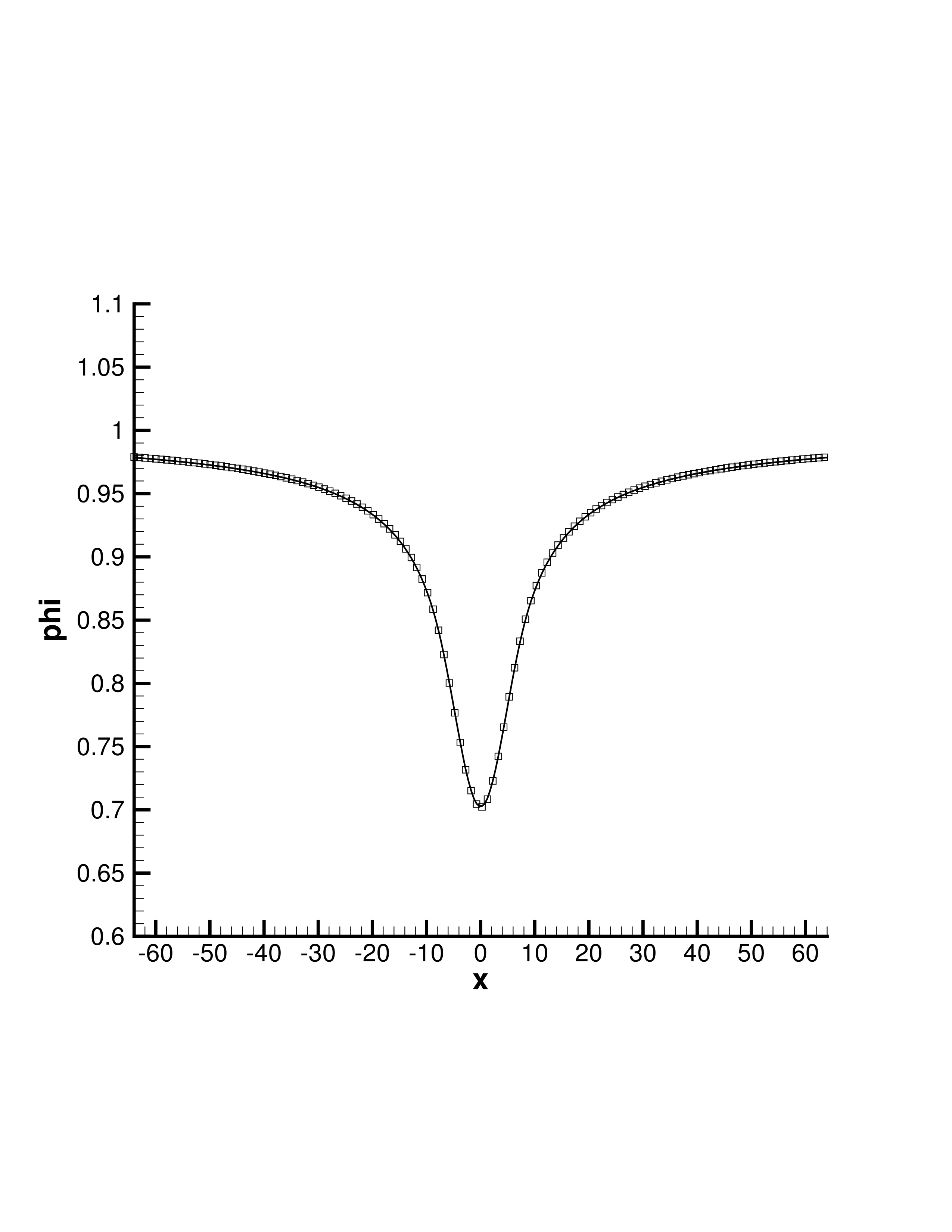} &  
		\includegraphics[width=0.45\textwidth]{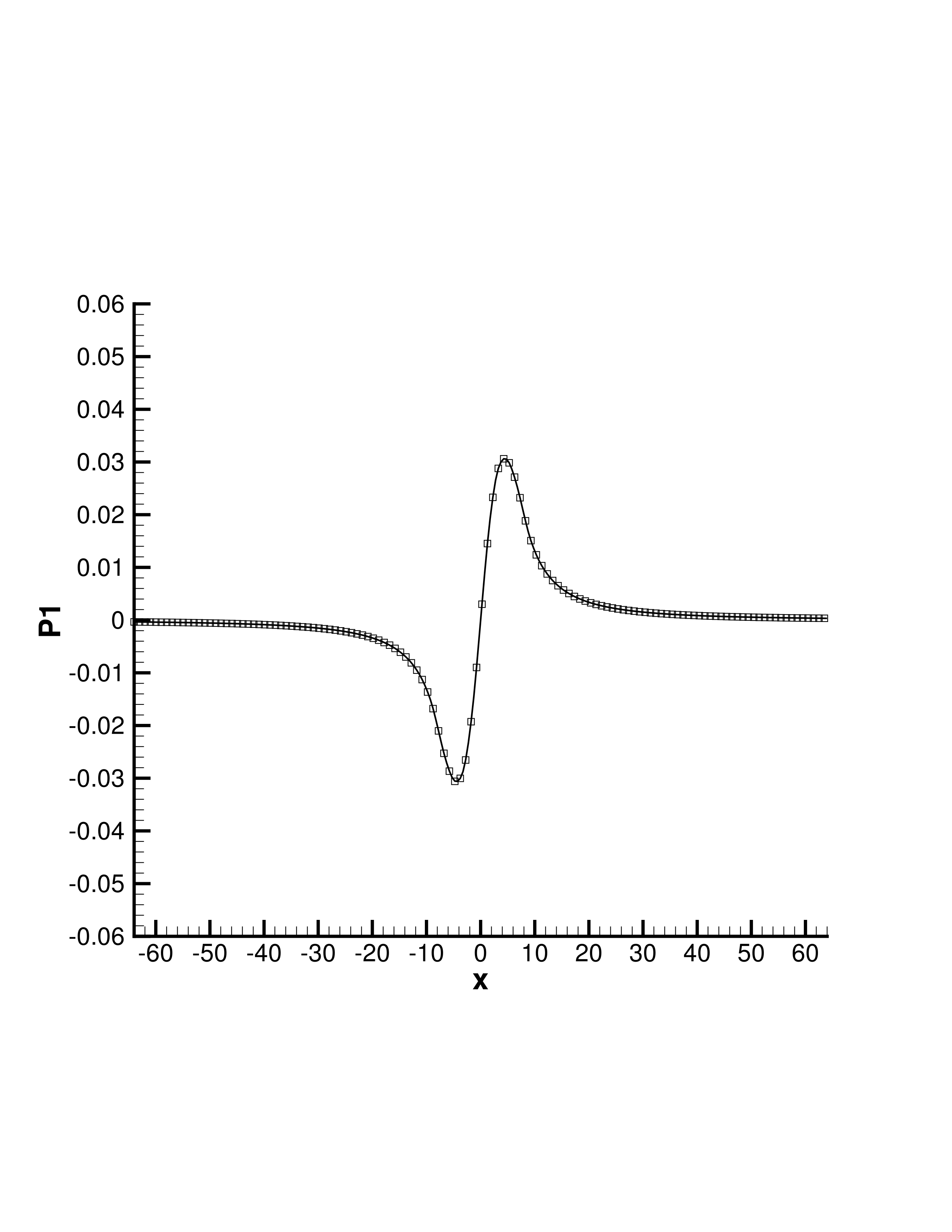}  \\ 
		\includegraphics[width=0.45\textwidth]{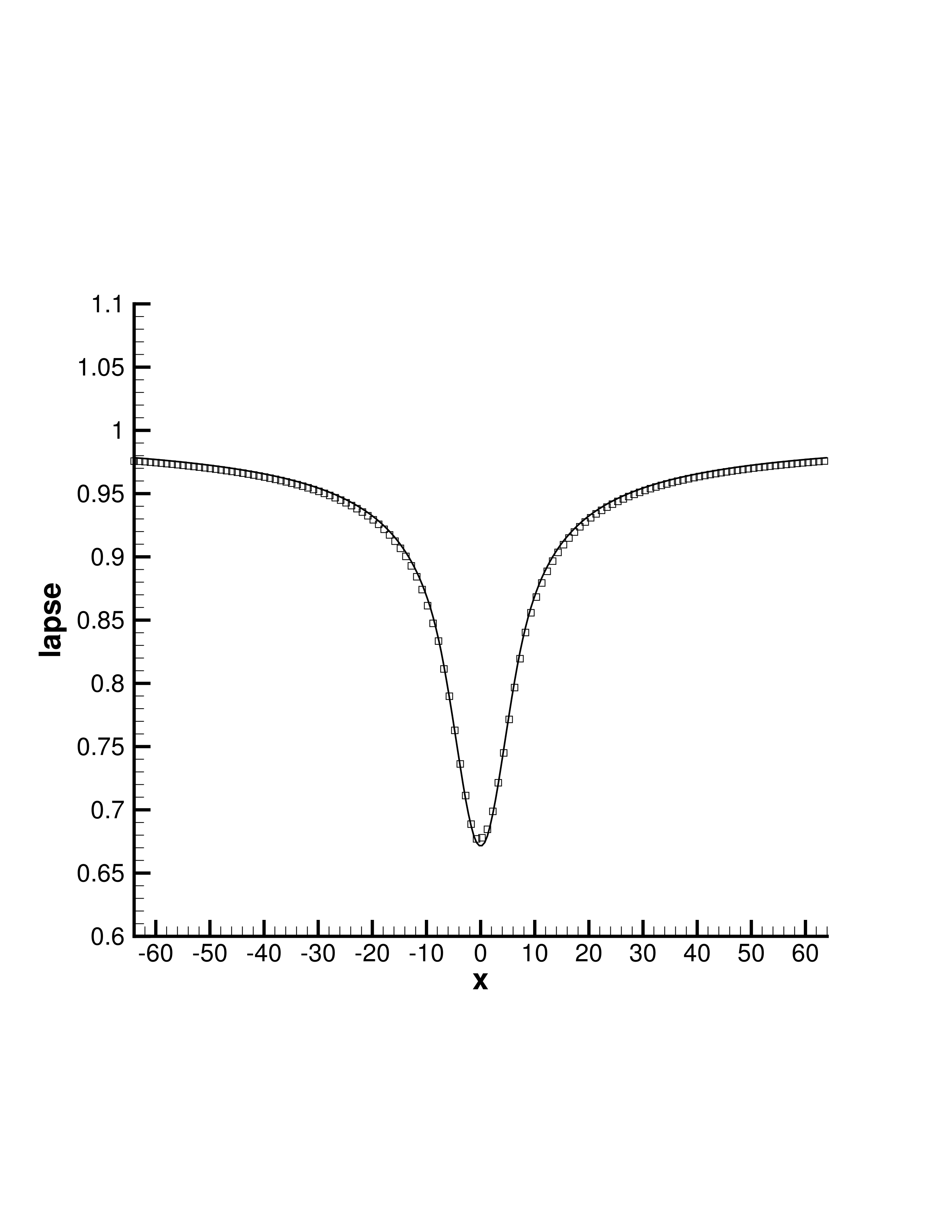} &  
	    \includegraphics[width=0.45\textwidth]{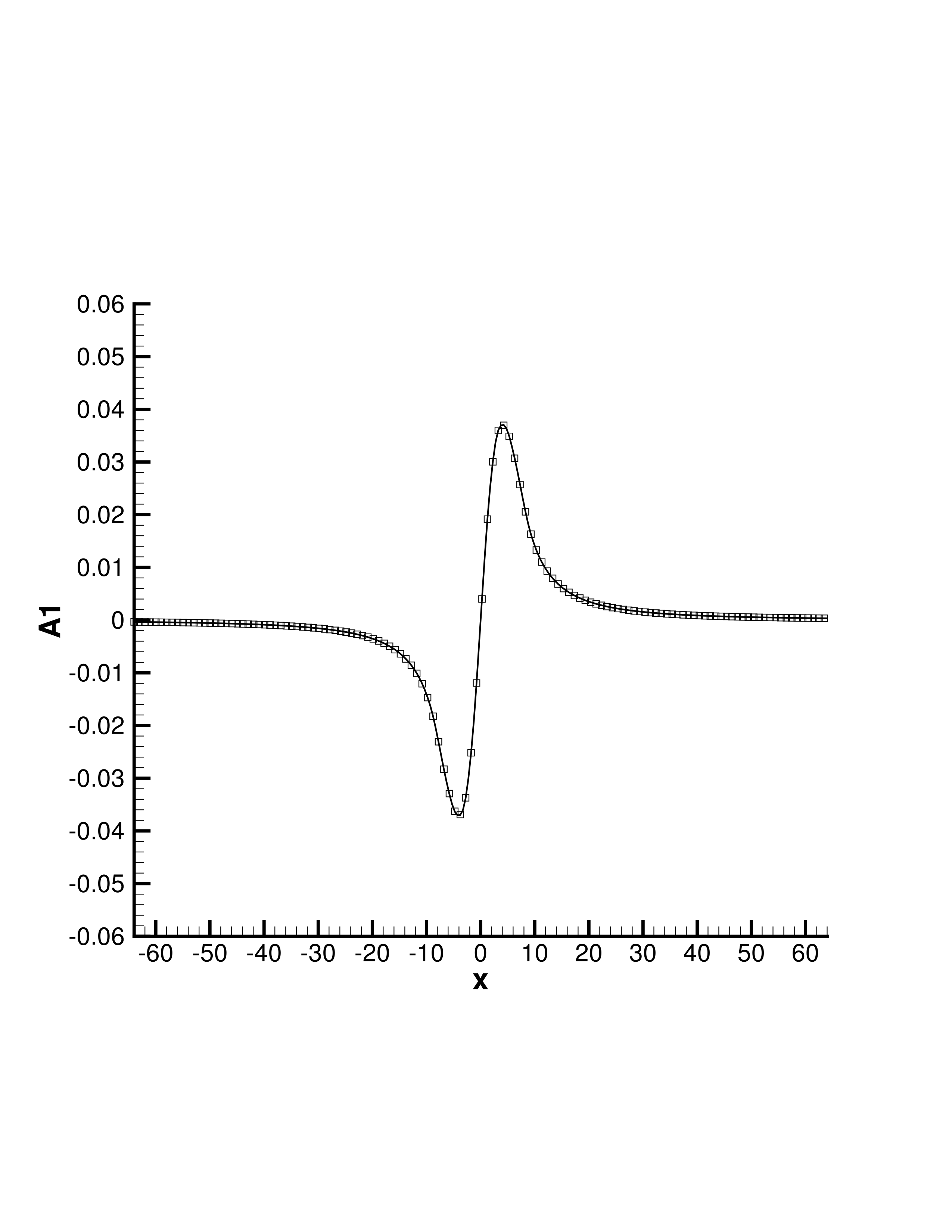}  \\ 
	\end{tabular} 
	\caption{Numerical simulation of a stable neutron star in anti-Cowling approximation at time $t=1000$. Comparison of the numerical solution (symbols) with the exact one (solid line) for the conformal factor $\phi$ (top left), for the lapse function $\alpha$ (bottom left) and for the associated auxiliary variables $P_1$ (top right) and $A_1$ (bottom right), which are subject to the curl constraints $\mathcal{P}_{lk}=0$ and $\mathcal{A}_{lk} = 0$, respectively.}
	\label{fig.metric}
\end{figure}

\begin{figure}
	\centering
		\includegraphics[width=0.7\textwidth]{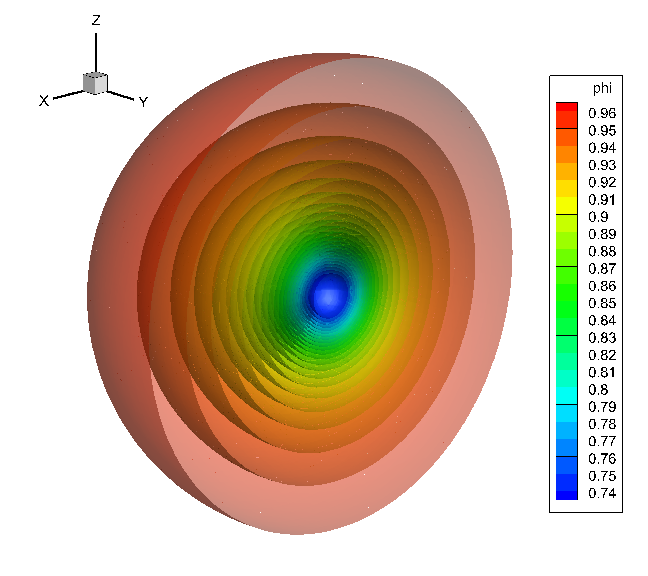}  
	\caption{Numerical simulation of a stable neutron star in anti-Cowling approximation at time $t=1000$ using the new hyperbolic GLM curl cleaning approach. Contour surfaces of the conformal factor $\phi$.}
	\label{fig.tov3d}
\end{figure}

The next test problem under consideration is the long time evolution of the spacetime generated by a stable 
neutron star (TOV star) in anti-Cowling approximation, i.e. 
we assume the matter quantities $\tau$ and $S_i$ to be \textit{stationary} in time and externally given 
by the Tolman-Oppenheimer-Volkoff (TOV) solution. For all details on the derivation of the radially 
symmetric TOV solution, which is used as initial condition for the present test case, see \cite{Tolman,Oppenheimer39b,Wald84,Carroll2003,bugner}. We run the test problem until a final time of $t=1000$ 
in the three-dimensional computational domain $\Omega=[-64,+64]^3$ using $63^3$ elements of polynomial 
approximation degree $N=3$.  
As gauge conditions we employ a frozen shift condition $\partial_t
\beta^i = 0$ by setting $s=0$ in the FO-CCZ4 system, together with the $1+\log$ slicing,
which corresponds to $g(\alpha)=2/\alpha$. The GLM cleaning speeds are set to $e=1.2$ for the  
cleaning of the nonlinear ADM constraints, while the cleaning speeds for the curl involutions
are chosen as $a_c^A = a_c^P = a_c^D = 0.1$, $a_d^A = a_d^D = a_d^P = 0.1$. The respective 
damping parameters are set to $\epsilon_{c,d}^A =  \epsilon_{c,d}^P = \epsilon_{c,d}^D= 5$. This means that only a very
small amount of GLM cleaning is used here. The reason for this choice is that the exact solution of the problem is smooth and \textit{stationary}, hence we expect only very small constraint violations to develop
due to the discretization errors of our fourth order ADER-DG scheme. 
The remaining constants in the FO-CCZ4 system are chosen as $\kappa_1=0.03$, $\kappa_2=\kappa_3=0$, $K_0=0$,
$c=0$ and $\eta=0$.   
The temporal evolution of the constraints without and with GLM cleaning are compared to each other in Figure
\ref{fig.constraints}, from which we can clearly conclude that even a very small amount of GLM cleaning 
is able to reduce the curl errors. This also means that despite the use of a very high order accurate scheme, the 
GLM cleaning is able to further reduce numerical errors in the constraints $\mathcal{A}_{lk}$, $\mathcal{P}_{lk}$ and
$\mathcal{D}_{lkij}$. In order to check the quality of our computational results at $t=1000$ in Figure \ref{fig.metric} 
we compare radial cuts through the numerical solution along the $x$ axis for the quantities $\phi$, $\alpha$, $P_1$ 
and $A_2$ with the exact solution, which is given by the initial condition. The agreement 
between numerical and exact solution is excellent for all quantities under consideration. Finally, in order
to visually check whether the numerical solution remains spherically symmetric, or not, we show iso contour
surfaces for the conformal factor $\phi$ at the final time $t=1000$ in Figure \ref{fig.tov3d}. From the
obtained results one can conclude that the solution remains clean and symmetric even after long integration 
times.

\subsection{Wavefield generated by two rotating masses}  
\label{sec.rotmass} 

\begin{figure}
	\centering
	\includegraphics[width=0.75\textwidth]{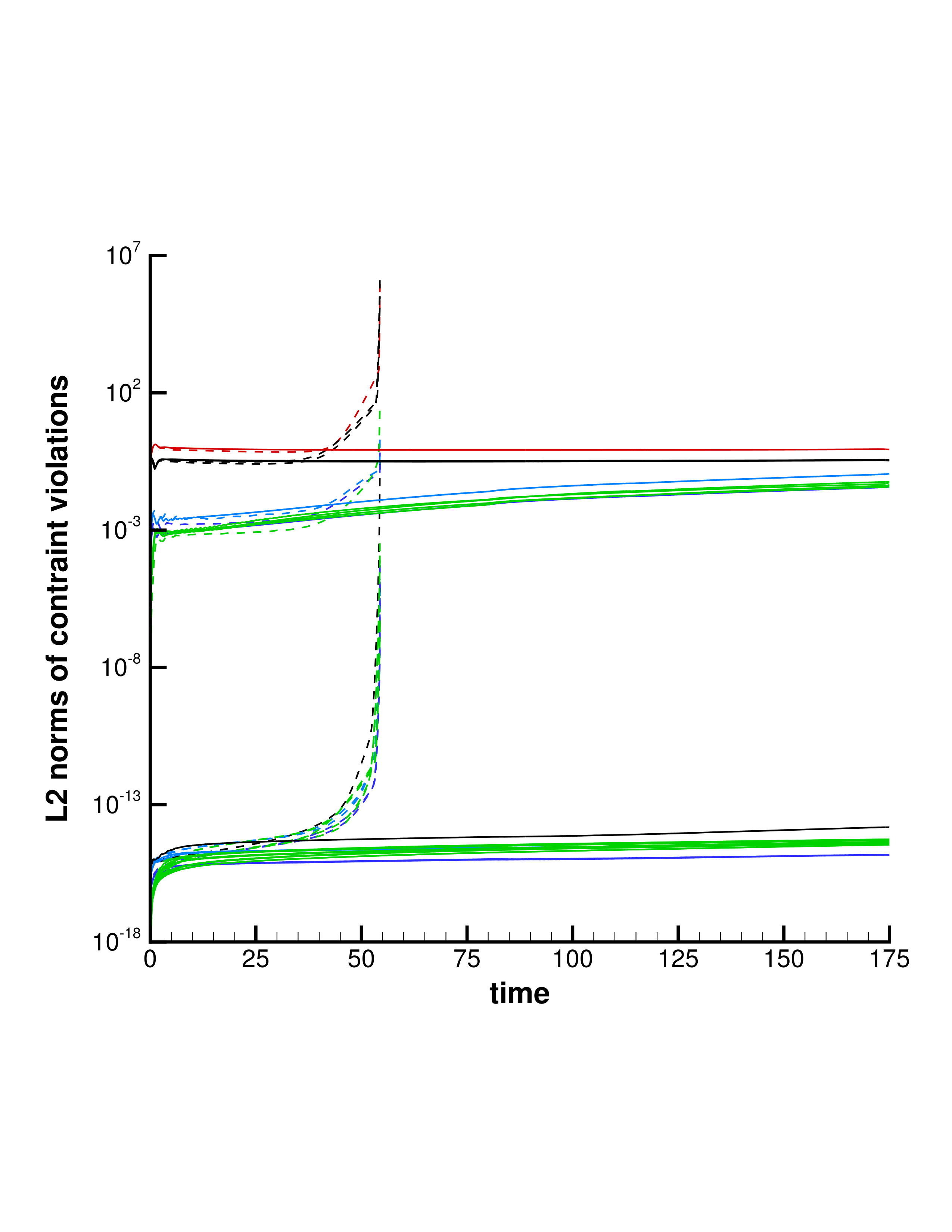}    	
	\caption{Constraint evolution for the rotating Gaussian density distributions. Dashed lines: standard FO-CCZ4 system without curl cleaning. Continuous lines: FO-CCZ4 with hyperbolic GLM curl cleaning. 
	The following constraints are monitored: Hamiltonian constraint $\mathcal{H}$ (red), momentum constraints $\mathcal{M}_i$ (black), curl constraints $\mathcal{A}_{lk}$ (light blue), $\mathcal{P}_{lk}$ (dark blue) and $\mathcal{D}_{lkij}$ (green).	
	}
	\label{fig.waveconst}
\end{figure}

\begin{figure}
	\centering
		\includegraphics[width=0.85\textwidth]{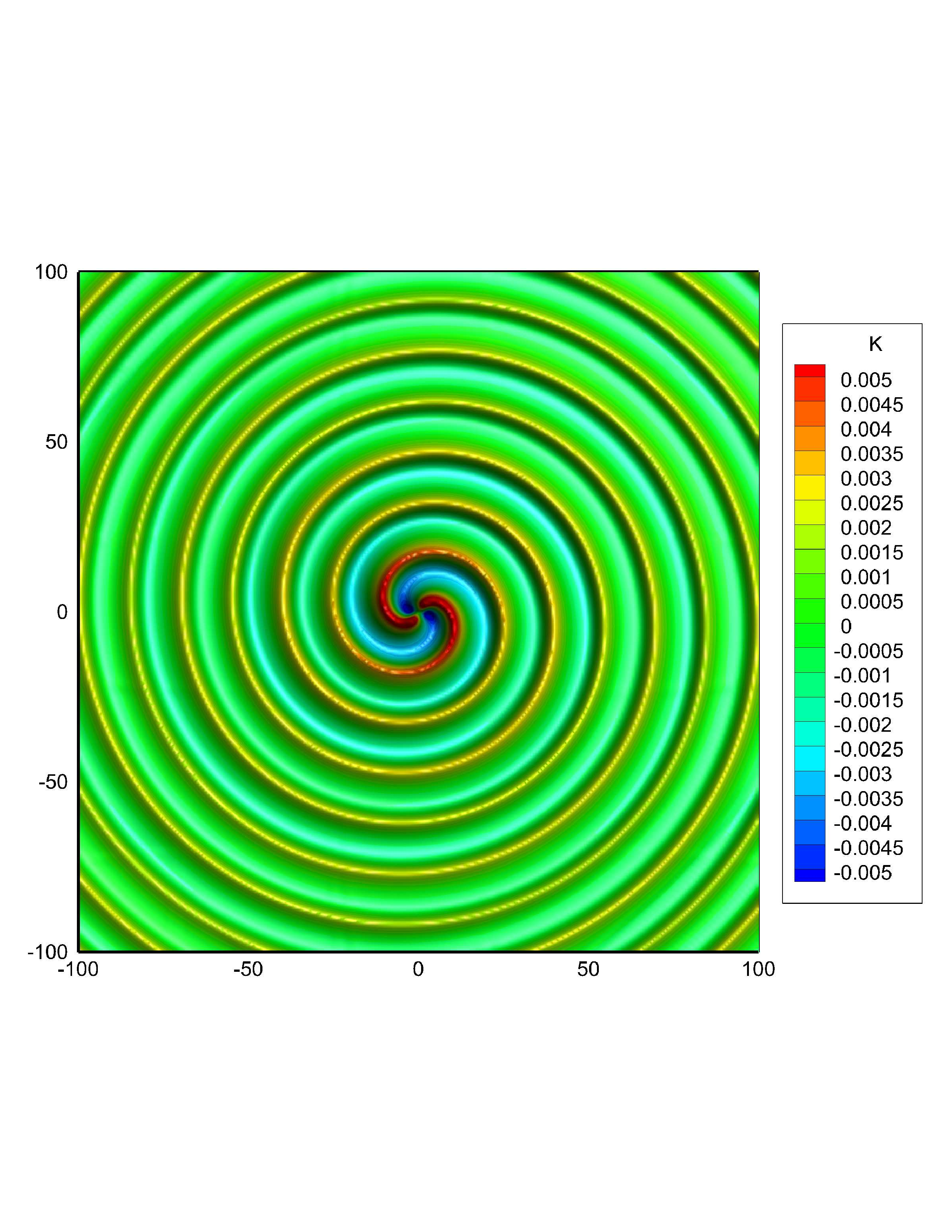}    	
	\caption{Contour colors of the trace of the extrinsic curvature ($K$) at time $t=175$ generated by two rotating Gaussian density distributions. }
	\label{fig.wavesK}
\end{figure}

This last test problem is only meant to be a showcase in order to demonstrate that the FO-CCZ4 system with the 
novel GLM curl cleaning proposed in this paper can also be used to simulate the propagation of waves in the space-time that are generated by the matter source terms. For this purpose, we start from a flat Minkowski spacetime and define an initial 
distribution of $\tau$ as 
\begin{equation}
\tau(\mathbf{x},0) = A_L \exp\left(-\frac{1}{2 \sigma_L^2}\left( \mathbf{x}-\mathbf{x}_L \right)^2 \right) + A_R \exp\left(-\frac{1}{2 \sigma_R^2}\left( \mathbf{x}-\mathbf{x}_R \right)^2 \right), 
\end{equation} 
which is then evolved in time via an artificially \textit{prescribed} motion 
\begin{equation}
  \partial_t \tau + v_k \, \partial_k \tau = 0 
\end{equation}
given by the background velocity field $\mathbf{v} = \boldsymbol{\Omega} \times \mathbf{x}$ for $\left\| \mathbf{x} \right\| < 5$ and $\mathbf{v}=0$ for $\left\| \mathbf{x} \right\| \geq 5$. Furthermore, 
in this test we set $S_i=0$ for all times and choose $A_L=A_R=5 \cdot 10^{-4}$, $\sigma_L = \sigma_R = 1$,
$\mathbf{x}_L=(-2,0,0)$, $\mathbf{x}_R=(2,0,0)$ and $\Omega=(0,0,0.2)$.  
The computational domain is  $\Omega=[-160,+160]^2 \times [-3.2,3.2]$ with periodic boundary conditions in 
$z$-direction, which would allow to solve the problem also in a purely two-dimensional setting. We employ 
a total of $320^2 \cdot 4$ elements with polynomial approximation degree $N=3$ in order to compute the
generated wavefield up to a final time of $t=175$. We run this test problem with two different setups. 
First, the standard FO-CCZ4 system \cite{ADERCCZ4} with the default choice $e=1$, $c=1$, $\kappa_i=0$ and 
without GLM curl cleaning is used. Then we run the same test problem again with GLM curl cleaning, 
using the cleaning speeds $e=2$, $a_c^A=a_c^P=a_c^D=a_d^A=a_d^P=a_d^D=1.5$, the damping 
parameters $\epsilon_c^A=\epsilon_c^P=\epsilon_c^D=\epsilon_d^A=\epsilon_d^P=\epsilon_d^D=1$ and
$c=0$, $\kappa_i=0$. 
In Figure \ref{fig.waveconst} we show the temporal evolution of the constraint violations obtained for
this test problem using a fourth order ADER-DG scheme ($N=3$). While the standard
FO-CCZ4 system quickly becomes highly unstable, the FO-CCZ4 system with GLM cleaning remains stable until 
the final time. In Figure \ref{fig.wavesK} we also show a snapshot of the generated wavefield by plotting 
the iso-contours of the extrinsic curvature at time $t=175$. The wavefield has a typical quadrupole-type
behaviour with its characteristic spirals, as expected from two rotating masses. 

\section{Conclusion} \label{sec.conclusion}

In this paper we have extended the GLM approach of Munz et al. \cite{MunzCleaning,Dedneretal}, which has originally only been developed for the divergence constraints in the Maxwell and MHD equations, to hyperbolic PDE systems with curl-type involutions. We have first presented the key ideas on a simple toy problem and have later extended the methodology to the first order reduction of the CCZ4 formulation of the Einstein field equations. The CCZ4 system is endowed with 11 natural curl involutions, namely one constraint 
on the variables $A_k$, another one on the variables $P_k$, three involutions for the field $B_k^{i}$ and six constraints for the variables $D_{kij}$. While the classical GLM cleaning for divergence-type constraints was 
based on a pair of div - grad operators, the new GLM approach proposed in this paper for the curl-type 
involutions makes use of pairs of curl - curl operators, which are therefore similar to Maxwell-type equations and thus are themselves endowed again with the usual divergence-free constraint that can again be taken into account by an additional GLM cleaning scalar. {For non-homogeneous curl involutions, even two additional cleaning scalars are needed. } 
We have shown several computational examples from computational general relativity in order to show that the proposed GLM cleaning can in some cases substantially improve the accuracy of the obtained computational results, even in the context of very high order accurate discontinuous Galerkin finite element schemes.

{In the near future, we plan to extend this approach also to the fully coupled Einstein-Euler system, i.e. where the matter source terms are given in a self-consistent manner by solving the general relativistic Euler or GRMHD equations, coupled with the FO-CCZ4 formulation of the Einstein field equations. }

 Further work will concern the application of the new GLM curl cleaning approach to other curl--constrained hyperbolic PDE systems, such as the novel hyperbolic surface tension model and the hyperbolic reformulation of the nonlinear Schr\"odinger equation of Gavrilyuk et al. \cite{HypSurfTension,Dhaouadi2018}, as well as to the compressible multi-phase model of Romenski et al. \cite{Rom1998,RomenskiTwoPhase2007,RomenskiTwoPhase2010}.  
 
 {Last but not least, the authors would also like to emphasize again that the proposed GLM curl cleaning approach is  meant to be only a first step into the direction of the design of curl-preserving numerical schemes. More work on the subject will be clearly necessary in the future, in particular concerning the development of high order accurate \textit{exactly} curl-preserving schemes on appropriately staggered meshes. The clear advantage of exactly curl preserving schemes would be not only that they require much less evolution variables, but also the fact that their structure-preserving property could be mathematically rigorously proven, while for GLM the structure is only asymptotically retrieved for large enough cleaning speeds. }
 
 {At this point we would like to emphasize again that the proposed GLM curl cleaning approach will only lead to a curl-preserving scheme in the asymptotic limit of infinitely fast cleaning speeds, which is obviously unfeasible in practice. It therefore remains up to the user to decide what level of curl errors are considered to be acceptable and which curl cleaning speeds are affordable for an explicit scheme due to the CFL stability condition. However, since the solutions of the FO-CCZ4 formulation of the Einstein field equations do not allow shock waves because all fields are linearly degenerate, see \cite{ADERCCZ4}, we believe that for the applications presented in this paper, and in combination with very high order discontinuous Galerkin finite element schemes, the proposed GLM curl cleaning approach is a very simple but appropriate technique for dealing with curl-constrained hyperbolic systems. }

\section*{Acknowledgments}
The research presented in this paper has been financed by the European Union's Horizon 2020 Research and  
Innovation Programme under the project \textit{ExaHyPE}, grant agreement number no. 671698 (call 
FETHPC-1-2014). 

M.D. also acknowledges the financial support received from the Italian Ministry of Education, University and Research (MIUR) 
in the frame of the Departments of Excellence Initiative 2018--2022 attributed to DICAM of the University of Trento 
(grant L. 232/2016) and in the frame of the PRIN 2017 project. M.D. has also received funding from the University of 
Trento via the  \textit{Strategic Initiative Modeling and Simulation}. 

E.G. has also been financed by a national mobility grant for young researchers in Italy, funded by GNCS-INdAM and acknowledges the support given by the University of Trento through the \textit{UniTN Starting Grant} initiative.

\section*{References}
\bibliographystyle{plain}
\bibliography{./GLMRot}

\begin{thebibliography}{10}

\bibitem{Alcubierre:2008}
M.~Alcubierre.
\newblock {\em Introduction to $3+1$ {N}umerical {R}elativity}.
\newblock Oxford University Press, Oxford, UK, 2008.

\bibitem{Alic:2009}
D.~{Alic}, C.~{Bona}, and C.~{Bona-Casas}.
\newblock {Towards a gauge-polyvalent numerical relativity code}.
\newblock {\em Phys. Rev. D}, 79(4):044026, 2009.

\bibitem{Alic:2012}
D.~{Alic}, C.~{Bona-Casas}, C.~{Bona}, L.~{Rezzolla}, and C.~{Palenzuela}.
\newblock {Conformal and covariant formulation of the Z4 system with
  constraint-violation damping}.
\newblock {\em Phys. Rev. D}, 85(6), 2012.

\bibitem{Arnowitt59}
R.~Arnowitt, S.~Deser, and C.~W. Misner.
\newblock Dynamical structure and definition of energy in general relativity.
\newblock {\em Phys. Rev.}, 116:1322, December 1959.

\bibitem{BalsaraAMR}
D.S. Balsara.
\newblock Divergence-free adaptive mesh refinement for magnetohydrodynamics.
\newblock {\em Journal of Computational Physics}, 174(2):614--648, 2001.

\bibitem{Balsara2004}
D.S. Balsara.
\newblock Second-order accurate schemes for magnetohydrodynamics with
  divergence-free reconstruction.
\newblock {\em The Astrophysical Journal Supplement Series}, 151:149--184,
  2004.

\bibitem{balsarahlle2d}
D.S. Balsara.
\newblock {Multidimensional HLLE Riemann solver: Application to Euler and
  magnetohydrodynamic flows}.
\newblock {\em Journal of Computational Physics}, 229:1970--1993, 2010.

\bibitem{balsarahllc2d}
D.S. Balsara.
\newblock {A two-dimensional HLLC Riemann solver for conservation laws:
  Application to Euler and magnetohydrodynamic flows}.
\newblock {\em Journal of Computational Physics}, 231:7476--7503, 2012.

\bibitem{MUSIC1}
D.S. Balsara.
\newblock {Multidimensional Riemann Problem with Self-Similar Internal
  Structure - Part I - Application to Hyperbolic Conservation Laws on
  Structured Meshes}.
\newblock {\em Journal of Computational Physics}, 277:163--200, 2014.

\bibitem{balsarahlle3d}
D.S. Balsara.
\newblock {Three dimensional HLL Riemann solver for conservation laws on
  structured meshes; Application to Euler and magnetohydrodynamic flows}.
\newblock {\em Journal of Computational Physics}, 295:1--23, 2015.

\bibitem{ADERdivB}
D.S. Balsara and M.~Dumbser.
\newblock {Divergence-free MHD on unstructured meshes using high order finite
  volume schemes based on multidimensional {R}iemann solvers}.
\newblock {\em Journal of Computational Physics}, 299:687--715, 2015.

\bibitem{MUSIC2}
D.S. Balsara and M.~Dumbser.
\newblock {Multidimensional Riemann Problem with Self-Similar Internal
  Structure - Part II - Application to Hyperbolic Conservation Laws on
  Unstructured Meshes}.
\newblock {\em Journal of Computational Physics}, 287:269--292, 2015.

\bibitem{BalsaraMultiDRS}
D.S. Balsara, M.~Dumbser, and R.~Abgrall.
\newblock {Multidimensional HLLC Riemann Solver for Unstructured Meshes - With
  Application to Euler and MHD Flows.}
\newblock {\em Journal of Computational Physics}, 261:172--208, 2014.

\bibitem{BalsaraGeodesic}
D.S. Balsara, V.~Florinski, S.~Garain, S.~Subramanian, and K.F. Gurski.
\newblock Efficient, divergence--free, high-order mhd on 3d spherical meshes
  with optimal geodesic meshing.
\newblock {\em Monthly Notices of the Royal Astronomical Society},
  487(1):1283--1314, 2019.

\bibitem{BalsaraCED}
D.S. Balsara, S.~Garain, A.~Taflove, and G.~Montecinos.
\newblock {Computational electrodynamics in material media with
  constraint-preservation, multidimensional Riemann solvers and sub-cell
  resolution -- Part II, higher order FVTD schemes}.
\newblock {\em Journal of Computational Physics}, 354:613--645, 2018.

\bibitem{BalsaraKaeppeli}
D.S. Balsara and R.~K\"appeli.
\newblock {von Neumann stability analysis of globally constraint-preserving
  DGTD and PNPM schemes for the Maxwell equations using multidimensional
  Riemann solvers}.
\newblock {\em Journal of Computational Physics}, 376:1108--1137, 2019.

\bibitem{BalsaraSpicer1999}
D.S. Balsara and D.~Spicer.
\newblock A staggered mesh algorithm using high order godunov fluxes to ensure
  solenoidal magnetic fields in magnetohydrodynamic simulations.
\newblock {\em Journal of Computational Physics}, 149:270--292, 1999.

\bibitem{Baumgarte99}
T.~W. {Baumgarte} and S.~L. {Shapiro}.
\newblock {Numerical integration of Einstein's field equations}.
\newblock {\em Phys. Rev. D}, 59(2):024007, 1999.

\bibitem{Baumgarte2010}
T.~W. {Baumgarte} and S.~L. {Shapiro}.
\newblock {\em {Numerical Relativity: Solving Einstein's Equations on the
  Computer}}.
\newblock Cambridge University Press, Cambridge, UK, 2010.

\bibitem{Bona:2003fj}
C.~{Bona}, T.~{Ledvinka}, C.~{Palenzuela}, and M.~{Z\'acek}.
\newblock {General-covariant evolution formalism for numerical relativity}.
\newblock {\em Phys. Rev. D}, 67(10):104005, May 2003.

\bibitem{Bona:2003qn}
C.~{Bona}, T.~{Ledvinka}, C.~{Palenzuela}, and M.~{Z\'acek}.
\newblock {Symmetry-breaking mechanism for the Z4 general-covariant evolution
  system}.
\newblock {\em Phys. Rev. D}, 69(6):064036, March 2004.

\bibitem{Bona95b}
C.~{Bona}, J.~{Mass{\'o}}, E.~{Seidel}, and J.~{Stela}.
\newblock {New Formalism for Numerical Relativity}.
\newblock {\em Phys. Rev. Lett.}, 75:600--603, 1995.

\bibitem{Bona:2004yp}
C.~{Bona} and C.~{Palenzuela}.
\newblock {Dynamical shift conditions for the Z4 and BSSN formalisms}.
\newblock {\em Phys. Rev. D}, 69(10):104003, May 2004.

\bibitem{Bona2009}
C.~Bona, C.~Palenzuela, and C.~Bona-Casas.
\newblock {\em Elements of Numerical Relativity and Relativistic Hydrodynamics:
  From Einstein's Equations to Astrophysical Simulations}.
\newblock Lecture Notes in Physics. Springer, Berlin Heidelberg, 2009.

\bibitem{Bona:2002fq}
C.~Bona, Ledvinka. T., and C.~Palenzuela.
\newblock A 3+1 covariant suite of numerical relativity evolution systems.
\newblock {\em Phys. Rev. D}, 66:084013, 2002.

\bibitem{ShashkovMultiMat3}
J.~Breil, T.~Harribey, P.~H. Maire, and M.J. Shashkov.
\newblock {A multi-material ReALE method with MOF interface reconstruction}.
\newblock {\em Computers and Fluids}, 83:115--125, 2013.

\bibitem{Brown2012}
J.~D. {Brown}, P.~{Diener}, S.~E. {Field}, J.~S. {Hesthaven}, F.~{Herrmann},
  A.~H. {Mrou{\'e}}, O.~{Sarbach}, E.~{Schnetter}, M.~{Tiglio}, and
  M.~{Wagman}.
\newblock {Numerical simulations with a first-order BSSN formulation of
  Einstein's field equations}.
\newblock {\em Phys. Rev. D}, 85(8):084004, 2012.

\bibitem{bugner}
M.~Bugner.
\newblock {\em Discontinuous Galerkin methods for general relativistic
  hydrodynamics}.
\newblock PhD thesis, Friedrich-Schiller-Universit\"at Jena, 2018.

\bibitem{Carroll2003}
S.~M. {Carroll}.
\newblock {\em {Spacetime and geometry. An introduction to general
  relativity}}.
\newblock 2003.

\bibitem{Dedneretal}
A.~Dedner, F.~Kemm, D.~Kr\"oner, C.~D. Munz, T.~Schnitzer, and M.~Wesenberg.
\newblock Hyperbolic divergence cleaning for the {MHD} equations.
\newblock {\em Journal of Computational Physics}, 175:645--673, 2002.

\bibitem{DelZanna2002}
L.~{Del Zanna} and N.~{Bucciantini}.
\newblock {An efficient shock-capturing central-type scheme for
  multidimensional relativistic flows. I. Hydrodynamics}.
\newblock {\em Astron. Astroph.}, 390:1177--1186, August 2002.

\bibitem{DeVore}
C.R. DeVore.
\newblock {Flux-corrected transport techniques for multidimensional
  compressible magnetohydrodynamics}.
\newblock {\em Journal of Computational Physics}, 92:142--160, 1991.

\bibitem{Dhaouadi2018}
F.~Dhaouadi, N.~Favrie, and S.~Gavrilyuk.
\newblock {Extended Lagrangian approach for the defocusing nonlinear
  Schr{\"{o}}dinger equation}.
\newblock {\em Studies in Applied Mathematics}, pages 1--20, 2018.

\bibitem{Axioms}
M.~Dumbser, F.~Fambri, M.~Tavelli, M.~Bader, and T.~Weinzierl.
\newblock {Efficient implementation of ADER discontinuous Galerkin schemes for
  a scalable hyperbolic PDE engine}.
\newblock {\em Axioms}, 7(3):63, 2018.

\bibitem{ADERCCZ4}
M.~Dumbser, F.~Guercilena, S.~K\"oppel, L.~Rezzolla, and O.~Zanotti.
\newblock {Conformal and covariant Z4 formulation of the Einstein equations:
  strongly hyperbolic first--order reduction and solution with discontinuous
  Galerkin schemes}.
\newblock {\em Physical Review D}, 97:084053, 2018.

\bibitem{GPRmodel}
M.~Dumbser, I.~Peshkov, E.~Romenski, and O.~Zanotti.
\newblock {High order ADER schemes for a unified first order hyperbolic
  formulation of continuum mechanics: Viscous heat-conducting fluids and
  elastic solids}.
\newblock {\em Journal of Computational Physics}, 314:824--862, 2016.

\bibitem{GPRmodelMHD}
M.~Dumbser, I.~Peshkov, E.~Romenski, and O.~Zanotti.
\newblock {H}igh order {ADER} schemes for a unified first order hyperbolic
  formulation of {N}ewtonian continuum mechanics coupled with electro-dynamics.
\newblock {\em Journal of Computational Physics}, 348:298--342, 2017.

\bibitem{AMR3DCL}
M.~Dumbser, O.~Zanotti, A.~Hidalgo, and D.S. Balsara.
\newblock {ADER-{WENO} Finite Volume Schemes with Space-Time Adaptive Mesh
  Refinement}.
\newblock {\em Journal of Computational Physics}, 248:257--286, 2013.

\bibitem{DGLimiter1}
M.~Dumbser, O.~Zanotti, R.~Loub{\`{e}}re, and S.~Diot.
\newblock {A posteriori subcell limiting of the discontinuous Galerkin finite
  element method for hyperbolic conservation laws}.
\newblock {\em Journal of Computational Physics}, 278:47--75, 2014.

\bibitem{Alcubierre:2003pc}
M.~Alcubierre et~al.
\newblock Toward standard testbeds for numerical relativity.
\newblock {\em Class. Quantum Grav.}, 21:589, 2004.

\bibitem{ArepoTN}
E.~Gaburro, W.~Boscheri, S.~Chiocchetti, C.~Klingenberg, V.~Springel, and
  M.~Dumbser.
\newblock {High order direct Arbitrary-Lagrangian-Eulerian schemes on moving
  Voronoi meshes with topology changes}.
\newblock {\em Journal of Computational Physics}, 2020.
\newblock submitted. https://arxiv.org/abs/1905.00967 ].

\bibitem{GardinerStone}
T.A. Gardiner and J.M. Stone.
\newblock {An unsplit Godunov method for ideal MHD via constrained transport}.
\newblock {\em Journal of Computational Physics}, 205:509--539, 2005.

\bibitem{GodunovRomenski72}
S.~K. Godunov and E.~I. Romenski.
\newblock Nonstationary equations of the nonlinear theory of elasticity in
  {Euler} coordinates.
\newblock {\em Journal of Applied Mechanics and Technical Physics},
  13:868--885, 1972.

\bibitem{Godunov:2003a}
S.~K. Godunov and E.~I. Romenski.
\newblock {\em {Elements of Continuum Mechanics and Conservation Laws}}.
\newblock Kluwer Academic/ Plenum Publishers, 2003.

\bibitem{Godunov1961}
S.K. Godunov.
\newblock An interesting class of quasilinear systems.
\newblock {\em Dokl. Akad. Nauk SSSR}, 139(3):521--523, 1961.

\bibitem{God1972MHD}
S.K. Godunov.
\newblock Symmetric form of the magnetohydrodynamic equation.
\newblock {\em Numerical Methods for Mechanics of Continuum Medium},
  3(1):26--34, 1972.

\bibitem{Gundlach:2005ta}
C.~Gundlach and J.M. Martin-Garcia.
\newblock Hyperbolicity of second-order in space systems of evolution
  equations.
\newblock {\em Class. Quantum Grav.}, 23:S387--S404, 2006.

\bibitem{Gundlach2005:constraint-damping}
Carsten Gundlach, Jose~M. Martin-Garcia, G.~Calabrese, and I.~Hinder.
\newblock Constraint damping in the {Z4} formulation and harmonic gauge.
\newblock {\em Class. Quantum Grav.}, 22:3767--3774, 2005.

\bibitem{HazraBalsara}
A.~Hazra, P.~Chandrashekar, and D.S. Balsara.
\newblock {Globally constraint-preserving FR/DG scheme for Maxwell's equations
  at all orders}.
\newblock {\em Journal of Computational Physics}, 394:298--328, 2019.

\bibitem{HymanShashkov1997}
J.M. Hyman and M.~Shashkov.
\newblock {Natural discretizations for the divergence, gradient, and curl on
  logically rectangular grids}.
\newblock {\em Computers and Mathematics with Applications}, 33:81--104, 1997.

\bibitem{JeltschTorrilhon2006}
R.~Jeltsch and M.~Torrilhon.
\newblock {On curl--preserving finite volume discretizations for shallow water
  equations}.
\newblock {\em BIT Numerical Mathematics}, 46:S35--S53, 2006.

\bibitem{MunzCleaning}
C.D. Munz, P.~Omnes, R.~Schneider, E.~Sonnendr\"ucker, and U.~Voss.
\newblock {Divergence Correction Techniques for Maxwell Solvers Based on a
  Hyperbolic Model}.
\newblock {\em Journal of Computational Physics}, 161:484--511, 2000.

\bibitem{Nakamura87}
T.~{Nakamura}, K.~{Oohara}, and Y.~{Kojima}.
\newblock {General Relativistic Collapse to Black Holes and Gravitational Waves
  from Black Holes}.
\newblock {\em Progress of Theoretical Physics Supplement}, 90:1--218, 1987.

\bibitem{Oppenheimer39b}
J.~R. Oppenheimer and G.~Volkoff.
\newblock On massive neutron cores.
\newblock {\em Phys. Rev.}, 55:374, 1939.

\bibitem{SHTC-GENERIC-CMAT}
I.~Peshkov, M.~Pavelka, E.~Romenski, and M.~Grmela.
\newblock {Continuum mechanics and thermodynamics in the Hamilton and the
  Godunov-type formulations}.
\newblock {\em Continuum Mechanics and Thermodynamics}, 30(6):1343--1378, 2018.

\bibitem{PeshRom2014}
I.~Peshkov and E.~Romenski.
\newblock A hyperbolic model for viscous {{N}ewtonian} flows.
\newblock {\em Continuum Mechanics and Thermodynamics}, 28:85--104, 2016.

\bibitem{GPRTorsion}
I.~Peshkov, E.~Romenski, and M.~Dumbser.
\newblock {Continuum mechanics with torsion}.
\newblock {\em Continuum Mechanics and Thermodynamics}, 2019.
\newblock in press. DOI: 10.1007/s00161-019-00770-6.

\bibitem{PowellMHD1}
K.G. Powell.
\newblock {An approximate Riemann solver for magnetohydrodynamics (that works
  in more than one dimension)}.
\newblock Technical Report ICASE-Report 94-24 (NASA CR-194902), NASA Langley
  Research Center, Hampton, VA, 1994.

\bibitem{PowellMHD2}
K.G. Powell, P.L. Roe, T.J. Linde, T.I. Gombosi, and D.L.~De Zeeuw.
\newblock A solution-adaptive upwind scheme for ideal magnetohydrodynamics.
\newblock {\em Journal of Computational Physics}, 154:284--309, 1999.

\bibitem{RomenskiTwoPhase2010}
E.~Romenski, D.~Drikakis, and E.F. Toro.
\newblock {Conservative models and numerical methods for compressible two-phase
  flow}.
\newblock {\em Journal of Scientific Computing}, 42:68--95, 2010.

\bibitem{RomenskiTwoPhase2007}
E.~Romenski, A.D. Resnyansky, and E.F. Toro.
\newblock {Conservative hyperbolic formulation for compressible two-phase flow
  with different phase pressures and temperatures}.
\newblock {\em Quarterly of Applied Mathematics}, 65:259--279, 2007.

\bibitem{Rom1998}
E.I. Romenski.
\newblock Hyperbolic systems of thermodynamically compatible conservation laws
  in continuum mechanics.
\newblock {\em Mathematical and computer modelling}, 28(10):115--130, 1998.

\bibitem{Schmidmayer2016}
K.~Schmidmayer, F.~Petitpas, E.~Daniel, N.~Favrie, and S.~Gavrilyuk.
\newblock {A model and numerical method for compressible flows with capillary
  effects}.
\newblock {\em Journal of Computational Physics}, 334:468--496, 2017.

\bibitem{HypSurfTension}
K.~Schmidmayer, F.~Petitpas, E.~Daniel, N.~Favrie, and S.~Gavrilyuk.
\newblock {Iterated upwind schemes for gas dynamics}.
\newblock {\em Journal of Computational Physics}, 334:468--496, 2017.

\bibitem{ShadabBalsara}
M.A. Shadab, D.S. Balsara, W.~Shyy, and K.~Xu.
\newblock Fifth order finite volume weno in general orthogonally -- curvilinear
  coordinates.
\newblock {\em Computers and Fluids}, 190:398--424, 2019.

\bibitem{Shibata95}
M.~{Shibata} and T.~{Nakamura}.
\newblock {Evolution of three-dimensional gravitational waves: Harmonic slicing
  case}.
\newblock {\em Phys. Rev. D}, 52:5428--5444, November 1995.

\bibitem{Springel}
V.~Springel.
\newblock {E pur si muove: Galilean-invariant cosmological hydrodynamical
  simulations on a moving mesh}.
\newblock {\em Monthly Notices of the Royal Astronomical Society (MNRAS)},
  401:791--851, 2010.

\bibitem{Tolman}
R.C. Tolman.
\newblock Static solutions of {Einstein's} field equations for spheres of
  fluid.
\newblock {\em Phys. Rev.}, 55:364--373, 1939.

\bibitem{Torrilhon2004}
M.~Torrilhon and M.~Fey.
\newblock {Constraint-preserving upwind methods for multidimensional advection
  equations}.
\newblock {\em SIAM Journal on Numerical Analysis}, 42:1694--1728, 2004.

\bibitem{Wald84}
Robert~M. Wald.
\newblock {\em General relativity}.
\newblock The University of Chicago Press, Chicago, 1984.

\bibitem{Yee66}
K.S. Yee.
\newblock {Numerical solution of initial voundary value problems involving
  Maxwell equation in isotropic media}.
\newblock {\em IEEE Trans. Antenna Propagation}, 14:302--307, 1966.

\bibitem{Zanotti2015b}
O.~{Zanotti}, F.~{Fambri}, and M.~{Dumbser}.
\newblock {Solving the relativistic magnetohydrodynamics equations with ADER
  discontinuous Galerkin methods, a posteriori subcell limiting and adaptive
  mesh refinement}.
\newblock {\em Mon. Not. R. Astron. Soc.}, 452:3010--3029, 2015.

\bibitem{Zanotti2015c}
O.~Zanotti, F.~Fambri, M.~Dumbser, and A.~Hidalgo.
\newblock Space-time adaptive {ADER} discontinuous {{G}alerkin} finite element
  schemes with a posteriori sub-cell finite volume limiting.
\newblock {\em Computers and Fluids}, 118:204 -- 224, 2015.

\end{thebibliography}

\clearpage

\end{document}